\documentclass[11pt,letterpaper]{amsart}
\usepackage[utf8]{inputenc}
\usepackage{float}
\usepackage{amsmath}
\usepackage{amsfonts}
\usepackage{amssymb}
\usepackage{amsthm}
\usepackage{graphicx}
\usepackage{lmodern}
\usepackage{tikz}
\usepackage{color}

\newtheorem{theorem}{Theorem}[section]

\newtheorem{lemma}[theorem]{Lemma}

\newtheorem{prop}[theorem]{Proposition}

\newtheorem{definition}[theorem]{Definition}

\newcommand{\C}{\mathbb{C}}
\newcommand{\R}{\mathbb{R}}
\newcommand{\RP}{\mathbb{RP}}
\newcommand{\Q}{\mathbb{Q}}

\newcommand{\N}{\mathbb{N}}

\newcommand{\Pa}{{\mathcal{P}}}
\newcommand{\s}{\mathbf{s}}
\newcommand{\bb}{\mathbf{b}}
\newcommand{\aalpha}{\boldsymbol{\alpha}}
\newcommand{\oomega}{\boldsymbol{\omega}}

\begin{document}

\title{Domains of Convergence for Polyhedral Packings}
\author[Ahmed et al.]{Nooria Ahmed, William Ball, Ellis Buckminster, Emilie Rivkin, Dylan Torrance, Jake Viscusi, Runze Wang, Ian Whitehead, S. Yang} 
\maketitle

\begin{abstract}
    Polyhedral circle packings are generalizations of the Apollonian packing. We develop the theory of the Apollonian group, Descartes quadratic form, and related objects for all polyhedral packings. We use these tools to determine the domain of absolute convergence of a generating function that can be associated to any polyhedral packing. This domain of convergence is the Tits cone for an infinite root system.
\end{abstract}

\section{Introduction}

Polyhedral packings are a class of plane circle packings which generalize the more well-known Apollonian circle packings, or Apollonian gaskets. In a standard Apollonian packing, four mutually tangent circles form what's known as a Descartes quadruple. Each Descartes quadruple has a dual quadruple, consisting of four circles which each pass through three points of tangency of the original circles. The Apollonian group is the group of M\"obius transformations generated by circular inversions through the dual circles. Its action on the original quadruple fills out the entire packing. We define the tangency graph of an Apollonian circle packing to be the graph with a vertex for each circle in the packing and an edge between vertices if the corresponding circles are tangent. Note that the tangency subgraph induced by a single Descartes quadruple will be the tetrahedral graph. By replacing this graph with the graphs of other polyhedra, we generate other polyhedral packings. For example, replacing the tetrahedron with the octahedron gives us a packing where the basic unit is a Descartes sextuple instead of a quadruple. In every polyhedral packing, the interiors of the circles are dense in the extended complex plane $\hat{\C}$, and the residual set or complement of these interiors is the limit set of a geometrically finite reflection group acting on the three-dimensional hyperbolic upper half-space. 

Examples of polyhedral packings have appeared in several articles. For example, the octahedral packing is introduced in \cite{MR2675919} and an asymptotic local-global property for curvatures in octahedral packings is proven in \cite{MR3781332}; the cubic packing is the $\Q[\sqrt{-2}]$ packing studied in \cite{MR3814328}. The general definition of polyhedral packings was introduced recently in \cite{MR3904690}, as part of a broader classification of crystallographic sphere packings. In \cite{MR3952499}, an asymptotic local-global property is proven at a level of generality that encompasses all superintegral polyhedral packings. We do not assume integrality in the present work, and we are only indirectly concerned with the question of which numbers can appear as curvatures in a packing. 

We have two main goals: first, to develop a theory of the geometric and algebraic Apollonian groups and the Descartes quadratic form for all polyhedral packings, generalizing the results of the influential papers \cite{MR1903421} and \cite{MR2173929}. Some of this theory is implicit in \cite{MR3904690} and other articles cited above, but we aim to make all the details explicit. Second, we introduce a generating series $Z(\s)$ which can be associated to any bounded polyhedral packing, generalizing \cite{Whitehead}. It is an exponential sum over all Descartes tuples that appear in the packing, and inherits a group of functional equations isomorphic to the Apollonian group. Analytic information about $Z(\s)$ can be used to study the density and other features of the set of tuples in the packing. We determine the domain of absolute convergence of this series, and visualize this domain as a subset of $\RP^3$--see Figure \ref{fig:domains}. It is the interior of the Tits cone of an infinite root system (sometimes a Kac-Moody root system) of complexity beyond affine and hyperbolic types. We show that the Tits cones of such root systems can encode the geometry of all polyhedral packing types. 

To give the formal construction of a polyhedral circle packing, begin with a 3-connected planar graph $\Pi$, which determines a polyhedron by Steinitz's theorem. By the Koebe-Andreev-Thurston theorem, this polyhedron can be embedded in 3-dimensional space with a midsphere $N$, i.e. a sphere tangent to all the edges of $\Pi$. There are two ways to use this embedding of $\Pi$ to generate a tuple of circles on the surface of $N$, which can be identified with $\hat{\C}$ by stereographic projection. First, from each vertex $i$ of $\Pi$, construct the unique cone tangent to $N$ with apex $i$. This cone intersects $N$ in a circle $c_i$. The collection of circles $c_i$ are pairwise tangent or disjoint, and can be oriented so that their interiors on $N$ are disjoint. They constitute the initial tuple for a polyhedral packing. The second construction begins with the collection of faces $j$ of $\Pi$. Each face $j$ intersects $N$ in a circle $d_j$. The circles $d_j$ are pairwise tangent or disjoint, and can be oriented so that their interiors on $N$ are disjoint. They are either disjoint from or orthogonal to the circles $c_i$. They constitute the dual tuple for the packing. 

\begin{figure}
\includegraphics[width=.72\textwidth]{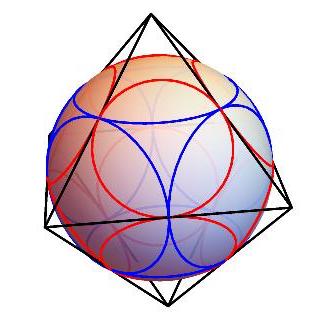}
\caption{An octahedron with a midsphere. The initial sextuple of circles for an octahedral packing is shown in blue on the sphere. The dual octuple of circles is shown in red. All figures in this article were created using Mathematica software \cite{Mathematica}.}
    \label{fig:midsphere}
\end{figure}

Figure \ref{fig:midsphere} shows an octhedron $\Pi$ with a midsphere $N$. The construction of the initial octahedral sextuple of circles and the dual cubic octuple are indicated. Note that the octahedron and the cube are dual to each other as polyhedra.

We record the precise tangency relations among circles $c_i$ and $d_j$, which are encoded by the combinatorics of $\Pi$. If vertices $i_1$, $i_2$ are connected by an edge of $\Pi$, then circles $c_{i_1}$ and $c_{i_2}$ are tangent at the point where this edge intersects $N$; otherwise, they are disjoint. If faces $j_1$, $j_2$ are adjacent along an edge of $\Pi$, then circles $d_{j_1}$ and $d_{j_2}$ are tangent at the point where this edge intersects $N$; otherwise, they are disjoint. If face $j$ contains vertex $i$, then circles $c_i$ and $d_j$ intersect orthogonally at the two points where edges through $i$ and $j$ intersect $N$; otherwise, $c_i$ and $d_j$ are disjoint. The tangency graph of the circles $c_i$ is the polyhedral graph $\Pi$; the tangency graph of the dual circles $d_j$ is the dual polyhedral graph. 

\begin{figure}
    \includegraphics[width=.72\textwidth]{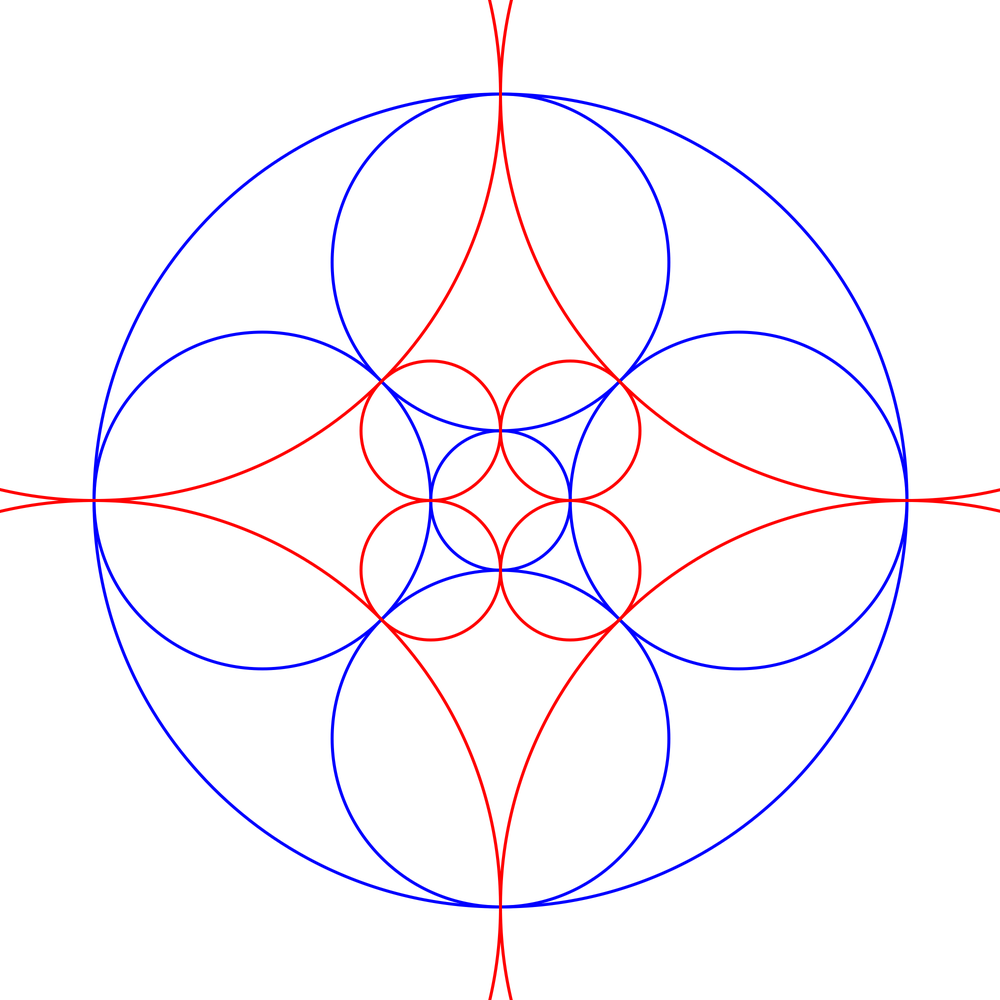}
    \caption{A Descartes sextuple of circles from an octahedral packing in blue and its dual, a Descartes octuple of circles from a cubic packing, in red.}
    \label{fig:duals}
\end{figure}

Figure \ref{fig:duals} gives an example of a sextuple of circles from an octahedral packing and its dual octuple of circles from a cubic packing in the plane. 

Inversion through a dual circle $d_j$ sends the original tuple of circles $c_i$ to another tuple of circles with the same tangency relations. Circles $c_i$ which intersect $d_j$ do so orthogonally, so they are preserved by this inversion. The group of M\"obius transformations generated by the inversions $\sigma_{j, \text{geom}}$ across dual circles $d_j$ is the geometric Apollonian group $W_{\text{geom}}$. This group has the presentation:
\begin{equation} W_{\text{geom}} = \langle \sigma_{1, \text{geom}}, \ldots \sigma_{n, \text{geom}} | \sigma_{j, \text{geom}}^2=I \rangle \end{equation}
We can see that there are no further relations because the circles $d_j$ have disjoint interiors. If $j_1, \ldots j_k$ is a list with consecutive terms distinct, and $z$ is a point not in the interior of any $d_j$, then $\sigma_{j_k, \text{geom}} \cdots \sigma_{j_1, \text{geom}}(z)$ lies in the interior of $d_{j_k}$, so $\sigma_{j_k, \text{geom}} \cdots \sigma_{j_1, \text{geom}}$ is not the identity. 

\begin{figure}
\includegraphics[width=.72\textwidth]{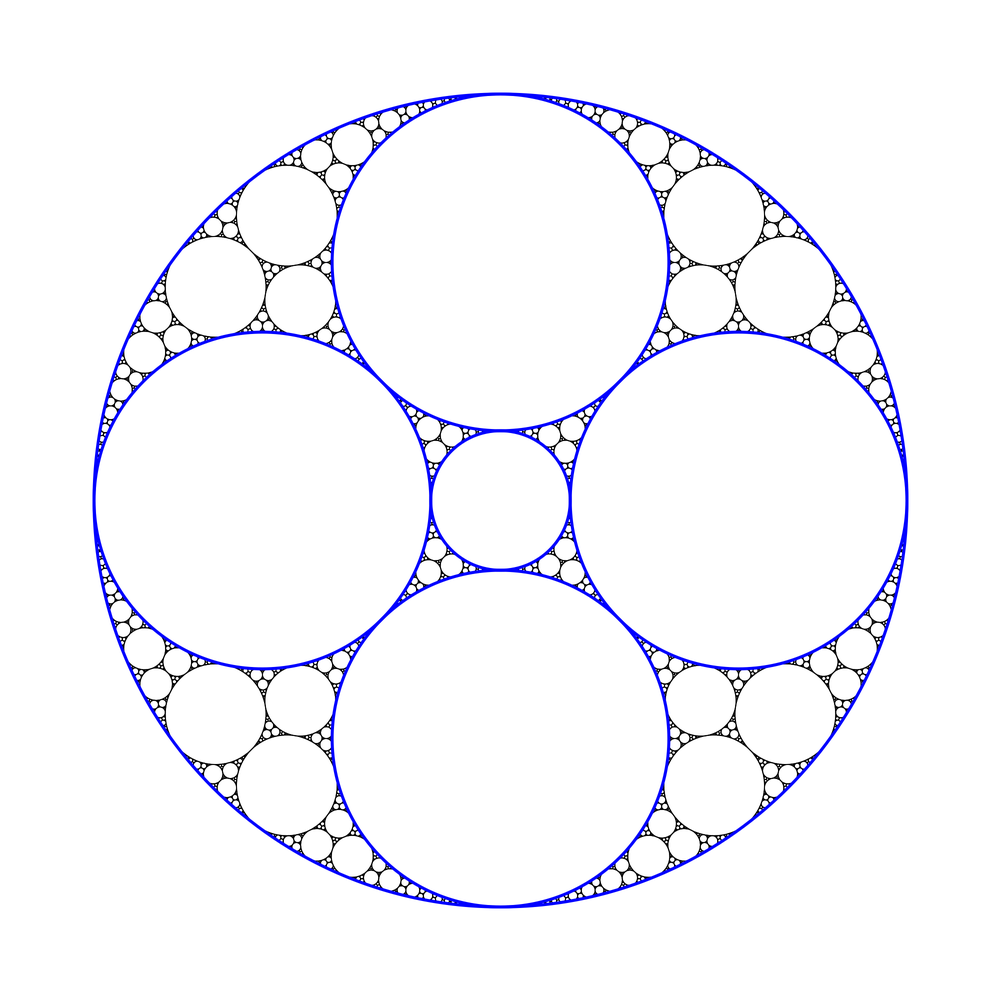}
\caption{An octahedral circle packing.}
\label{fig:packing}
\end{figure}

The polyhedral circle packing $\Pa$ is the orbit of $W_{\text{geom}}$ on the initial tuple of circles $c_i$. Figure \ref{fig:packing} shows an octahedral packing generated from the initial circles and dual circles shown in Figure \ref{fig:duals}. Note that we could just as easily have generated a cubic packing by reflecting the dual circles across the initial tuple. 

To see how the polyhedral packings generalize Apollonian packings, note that in the standard Apollonian packing, the basic unit is a Descartes quadruple of circles, whose tangency graph is a tetrahedron. The four dual circles whose inversions generate the standard Apollonian group can be identified with the faces of the tetrahedron. In this case, the tetrahedron is its own dual polyhedron. 

We now introduce the generating series which will be one of our main objects of study. In \cite{Whitehead}, the series $Z(\mathbf{s})$ which can be seen as a generating function for the tetrahedral Apollonian packing, is introduced. We define a generalized version for all polyhedral packings:

\begin{definition}\label{Z(s)}
For a fixed bounded polyhedral packing $\Pa$, define $Z \colon \C^m \to \C$ as
\begin{equation} Z(\s) = \sum_{\bb\in \Pa} e^{-\bb\cdot \s} \end{equation}
where $\mathbf{s} = (s_1, \ldots s_m) \in \C^m$ and the sum runs over all Descartes $m$-tuples of curvatures $\bb = (b_1, \ldots b_m) \in \Pa$, with multiplicity if a tuple appears more than once in $\Pa$.
\end{definition}

The convergence of $Z(\s)$ gives information about the growth of the set of $m$-tuples in $\Pa$. Following \cite{Whitehead}, a Mellin transform of $Z$ can be used to produce the L-series associated to $\Pa$, whose abscissa of convergence is the Hausdorff dimension of $\Pa$. Moreover, $Z(\s)$ encapsulates much of the structure of $\Pa$--we will see that it has an infinite group of functional equations isomorphic to $W_{\text{geom}}$. This symmetry group is crucial in determining the domain of absolute convergence. 

In Section 2, we give analogues of the Descartes quadratic form and algebraic Apollonian group for all polyhedral packings. Descartes found that quadruples of curvatures $\bb=(b_1, b_2, b_3, b_4)$ of mutually tangent circles in a tetrahedral configuration satisfy the quadratic equation 
\begin{equation}2(b_1^2+b_2^2+b_3^2+b_4^2)-(b_1+b_2+b_3+b_4)^2=0 \end{equation}
For the tetrahedral packing, the algebraic Apollonian group $W$ is isomorphic to the geometric Apollonian group $W_{\text{geom}}$, but  realized as a group of $4 \times 4$ matrices which preserve the Descartes quadratic form. It permutes the collection of Descartes quadruples in a packing. We follow the methods of \cite{MR1903421} and \cite{MR2173929} to extend these ideas. We find that all $m$-tuples of curvatures $\bb$ in a polyhedral packing $\Pa$ lie in a $4$-dimensional subspace of $\R^m$, and satisfy a homogeneous quadratic equation defined on this subspace. Theorem \ref{descartes} is our generalization of the Descartes circle theorem. Further, there is an algebraic Apollonian group $W$ isomorphic to $W_{\text{geom}}$ which permutes the collection of tuples in $\Pa$ and preserves the quadratic form. We give detailed constructions of all these objects from the data of the initial tuple and dual tuple for $\Pa$. In the cases of the octahedral and cubic packings, we recover results of Guettler and Mallows \cite{MR2675919} and Stange \cite{MR3814328} respectively. 

In Sections 3 and 4, we turn our attention to the domain of convergence of $Z(\s)$. From the fact that $m$-tuples $\bb \in \Pa$ lie in a $4$-dimensional subspace, we can consider $Z$ as a function on the dual $4$-dimensional quotient of $C^m$. Furthermore, since $Z(\s)$ is an exponential series, its domain of convergence is a convex cone, which means that it can be visualized in $3$-dimensional projective space. Since the algebraic Apollonian group $W$ acts strictly transitively on tuples of circles in $\Pa$, $Z(\s)$ can be rewritten as a symmetric function or the sum over an orbit of this group. This implies that it is invariant under the action of the transpose group $W^T$ on $\s$. 

In Section 3, we use lower bounds on the growth of tuples in $\Pa$ to establish an initial domain of convergence $A_0$ for $Z(\s)$. In $\RP^3$, $A_0$ is a polyhedron of the same type that defines $\Pa$. Theorem \ref{finaldomain} uses the $W^T$-invariance to describe the full domain of absolute convergence $A$: 
\begin{equation} A=\bigcup_{w\in W} w^T A_0 \end{equation}
In Section 4, our subject is the geometry of $A$. Theorem \ref{domaingeometry} states that, in $\RP^3$, $A$ is the union of a ball with infinitely many solid cones, tangent to the ball, whose bases form a packing of the same type as $\Pa$ on the surface of the ball. For illustrations of $A$ and its construction as the $W^T$-orbit of $A_0$, see Figure \ref{fig:domains}. 

This article contributes to the study of general polyhedral packings initiated in \cite{MR3904690}. It is natural to ask whether our results extend to the broader class of crystallographic sphere packings introduced in that article. We conjecture that all our results here are true for crystallographic packings in full generality. Some of the underlying definitions need to be modified. One difficulty is that there is no canonical definition of a tuple or dual tuple of spheres in a crystallographic packing. It becomes necessary to make an arbitrary choice of sufficiently many spheres in $\Pa$ as the tuple (more precisely, a collection of spheres whose augmented center-curvature coordinates span $\R^{D+2}$, $D$ being the dimension of the packing). Any such choice will give rise to isomorphic versions of the root space and weight space, the bilinear forms on these spaces, and the algebraic Apollonian group. The definition of $Z(\s)$, even the number of variables, will depend on the choice of tuple, but the final domain of convergence, as a subset of the weight space, will be uniquely determined up to a linear change of variables. 

We also aim to contribute to a broader research program with many open questions: classifying the projective geometry of infinite root systems, including limit roots, limit weights, Tits buildings, imaginary cones, and Tits cones. These objects have been explored in \cite{MR4036724}, \cite{MR3563253}, \cite{MR3176144}, \cite{MR4104416}, and elsewhere. Here we must make some comments on how root system terminology is used in this article. Our root systems satisfy the following very general definition, used in \cite{MR3563253} and \cite{MR3176144}. In a real vector space $V$ equipped with a symmetric bilinear form $\tilde{G}$, a finite set of simple roots $\aalpha_1, \ldots \aalpha_n$ must satisfy the following properties:
\begin{itemize}
    \item The $\aalpha_j$ are positively independent, i.e. if $\sum a_j \aalpha_j = 0$ with all $a_j \geq 0$, then all $a_j=0$.
    \item For $j_1 \neq j_2$, $\aalpha_{j_1}^T \tilde{G} \aalpha_{j_2} \in (-\infty, -1] \cup \lbrace -\cos\left( \frac{\pi}{k} \right) | k=2, 3, 4, \ldots \rbrace$.
    \item For all $j$, $\aalpha_{j}^T \tilde{G} \aalpha_{j}=1$. 
\end{itemize}
Then the Weyl group $W$ is the group generated by reflections with respect to the simple roots; the orbit of $W$ on the simple roots is the full collection of roots. We also employ a definition of weights adapted to this scenario. The fundamental weights are elements of the dual space $V^*$ to $V$. They are not a dual basis to $\lbrace \aalpha_j \rbrace$ (since these may not be linearly independent); rather, the fundamental weights are vertices of the dominant cone, whose faces correspond to simple roots. The fundamental weights inherit an action of $W$, and their orbit is the full collection of weights. 

Let us give a dictionary between root system and polyhedral packing terminologies used in this article. See Section 2 for the relevant definitions. In our setup, the root space is the space spanned by Descartes tuples in a packing; the bilinear form on this space is the Descartes quadratic form. The weights are in $1:1$ correspondence with the circles in a packing, with fundamental weights corresponding to circles in an initial tuple. The roots correspond to dual circles, with simple roots corresponding to circles in the dual tuple. The Weyl group is precisely the Apollonian group for the packing. The dominant cone for the root system is our initial cone of convergence $A_0$, and the interior of the Tits cone is the full cone of convergence $A$.

Our root systems do not all satisfy the integrality conditions needed to be Kac-Moody root systems, but some do. In particular, if the Gram matrix $G$ of the dual polyhedron can be obtained from an integral Cartan matrix, i.e. if $2DGD^{-1}$ has integer entries for some diagonal matrix $D$, then the root system is Kac-Moody. This will hold for integral polyhedral packings, like the cubic and octahedral packings. 

In \cite{MR679972}, Maxwell gives a classification of sphere packings by root systems of ``level 2.'' These are root systems with Coxeter diagrams where, if any one vertex is removed, the result is a finite, affine or hyperbolic diagram; if two vertices are removed, the result is finite or affine. Chen and Labb\'{e} show that, in these root systems, the residual set of the packing is realized as the set of limit weights or limit roots \cite{MR3303040}. The limit weights lie in the same space as the projective Tits cone that we study, and the two structures are closely related. Maxwell tabulates the root systems which give rise to packings; Chen and Labb\'{e} confirm and extend these tables by computer calculation. We approach their question from the opposite perspective: instead of asking which root systems lead to packings, we ask which packings can appear in the Tits cone of some root system. Our Theorem \ref{domaingeometry} implies that every polyhedral packing appears in the Tits cone of a root system, sometimes a Kac-Moody root system. 

This suggests that the world of limit weights and Tits cones for root systems encompasses the whole world of circle and sphere packings. But the world of limit weights and Tits cones is much larger. There are infinitely many level 2 root systems in ranks 3 and 4 (polyhedral packings have rank 4), but finitely many in ranks 5 and above, peaking at rank 11. Above level 2, very little is known about the structure of the limit weights and Tits cone. The signature of the quadratic form is no longer necessarily $(n,1)$, so the techniques of hyperbolic geometry may not apply. The limit weights and Tits cones of higher level root systems are a vast class of unknown fractal sets, generalizing circle and sphere packings. We hope to study these objects more in future work.

\textbf{Acknowledgements:} We are grateful to Arthur Baragar, Lisa Carbone, Christophe Hohlweg, Cathy Hsu, Edna Jones, Alex Kontorovich, and Kate Stange, for interesting discussions related to this project.

\section{Analogues of the Descartes Quadratic Form and Apollonian Group for Polyhedral Packings} 
In this section we give analogues of the Descartes quadratic form and the Apollonian group for all polyhedral packings. In doing so, we introduce notation and formalism that will be used throughout this article. 

To describe generalized circles and packings in the extended complex plane we use the augmented curvature-center coordinates described in \cite{MR1903421}. For a fuller exploration of this coordinate system, see \cite{Kocik2}.  Each circle $c$ is represented by a vector of the form $\mathbf{c}=(\tilde{b},b,h_1,h_2)^T \in \R^4$, where $b$ is the curvature, $\tilde{b}$ is the curvature after inversion through the unit circle and $h_1,h_2$ are $b$ times the $x$ and $y$ coordinates of the center, respectively. The circle is oriented with normal vector pointing inward if $b>0$ and outward if $b<0$. If the curvature $b$ is $0$,  then the circle is a line, and $h_1$ and $h_2$ are the $x$ and $y$ coordinates of the unit normal vector. In this case, the direction of $(h_1, h_2)$ gives the orientation. Each generalized circle divides $\hat{\C}$ into two disjoint regions; we say that the interior of the circle is the region that the normal vector points toward. 

We have the bilinear form 
\begin{equation} P=\begin{pmatrix} 0 & -\tfrac{1}{2} & 0 & 0 \\ \tfrac{1}{2} & 0 & 0 & 0 \\ 0 & 0 & 1 & 0 \\ 0 & 0 & 0 &1 \end{pmatrix} \end{equation} 
of signature $(3, 1)$ on the space of circles. This form has a nice geometric interpretation:
\begin{prop}[\cite{Kocik2}, Prop. 2.4]\label{bilinearform}
We have $\mathbf{c}^T P \mathbf{c}=1$ for all vectors $\mathbf{c}$ representing oriented generalized circles in $\hat{\C}$. Furthermore, if $\mathbf{c}_1, \mathbf{c}_2$ represent two distinct oriented generalized circles in $\hat{\C}$, then 

\noindent\begin{tabular}{@{}ll}
$\mathbf{c}_1^T P \mathbf{c}_2<-1$ & if $\mathbf{c}_1, \mathbf{c}_2$ are disjoint, neither interior contains the other \\
$\mathbf{c}_1^T P \mathbf{c}_2=-1$ & if $\mathbf{c}_1, \mathbf{c}_2$ are tangent, neither interior contains the other \\
$\mathbf{c}_1^T P \mathbf{c}_2=\cos(\theta)$ & if $\mathbf{c}_1, \mathbf{c}_2$ intersect, where $\theta$ is the angle between normal \\ \, & vectors at a point of intersection \\
$\mathbf{c}_1^T P \mathbf{c}_2=0$ & if $\mathbf{c}_1, \mathbf{c}_2$ intersect orthogonally \\
$\mathbf{c}_1^T P \mathbf{c}_2=1$ & if $\mathbf{c}_1, \mathbf{c}_2$ are tangent, one interior contains the other \\
$\mathbf{c}_1^T P \mathbf{c}_2>1$ & if $\mathbf{c}_1, \mathbf{c}_2$ are disjoint, one interior contains the other \\
\end{tabular}
\end{prop}
In addition, the augmented center-curvature coordinate system interacts well with M\"{o}bius transformations:
\begin{prop}[\cite{MR2173929}, Thm. 2.8]\label{mobius}
The action of M\"{o}bius transformations on generalized circles in $\hat{\C}$ is linear in the augmented curvature-center coordinate system, and preserves the bilinear form $P$. The resulting homomorphism from the group of M\"{o}bius transformations to $O_P$ is injective.
\end{prop}

Let $\mathbf{c}_1, \ldots \mathbf{c}_m$ be a polyhedral configuration of circles and let $\mathbf{d}_1, \ldots \mathbf{d}_n$ be the dual circles. The values of the bilinear form $P$ on $\mathbf{c}_i$ and $\mathbf{d}_j$ give all the information we will need in our calculation. By Proposition \ref{bilinearform}, we find that: 
\begin{align}
&\begin{array}{ll} \mathbf{c}_i^T P \mathbf{c}_j =1 & \text{ if } i=j \\ \mathbf{c}_i^T P \mathbf{c}_j =-1 & \text{ if } i, j \text{ are adjacent} \\ \mathbf{c}_i^T P \mathbf{c}_j <-1 & \text{ otherwise} \end{array} \\
&\begin{array}{ll} \mathbf{d}_i^T P \mathbf{d}_j =1 & \text{ if } i=j \\ \mathbf{d}_i^T P \mathbf{d}_j =-1 & \text{ if } i, j \text{ are adjacent} \\ \mathbf{d}_i^T P \mathbf{d}_j <-1 & \text{ otherwise} \end{array} \label{config1} \\
&\begin{array}{ll} \mathbf{c}_i^T P \mathbf{d}_j =0 & \text{ if face j contains vertex i} \\ \mathbf{c}_i^T P \mathbf{d}_j <-1 & \text{ otherwise} \end{array} \label{config2}
\end{align}

By Proposition \ref{mobius}, the group of Mobius transformations generated by reflections $\sigma_{1,\text{geom}}, \ldots \sigma_{n,\text{geom}}$ across the dual circles can be identified with its image in $O_P$. We call this group the geometric Apollonian group $W_{\text{geom}}$. For each $j$, 
\begin{equation}\sigma_{j,\text{geom}}(\mathbf{v})=\mathbf{v}- 2(\mathbf{d}_j^TP\mathbf{v}) \mathbf{d}_j \end{equation}
As shown in the previous section, the only relations among the generators are that each $\sigma_{j,\text{geom}}^2=I$. 

A polyhedral circle configuration can be identified with a $4 \times m$ matrix: 
\begin{equation} C= \left( \begin{array}{cccc} | & | & \, & | \\ \mathbf{c}_1 & \mathbf{c}_2 & \cdots & \mathbf{c}_m \\ | & | & \, & | \end{array}\right) \end{equation}
We will need the following lemma:
\begin{lemma} \label{rank4}
$C$ has rank $4$. 
\end{lemma} 
\begin{proof}
By applying a M\"{o}bius transformation, we may send two tangent circles in the configuration to the lines $\Im (z)=0$ and $\Im (z) =1$ in $\hat{\C}$, and one of their common dual circles to the line $\Re (z)=0$. The coordinates of the two tangent circles become $(0,0,0,-1)^T$ and $(2, 0, 0, 1)^T$. Any additional circle orthogonal to the dual circle will have $h_1=0$ but $b \neq 0$, so it will be linearly independent from these two. Any additional circle not orthogonal to the dual circle will have $h_1 \neq 0$ so it will be linearly independent from these three.
\end{proof}

A packing can be described as the orbit of $W_{\text{geom}}$ on an initial tuple $C$. By the lemma, $C$ admits a right inverse, so $w \in W_{\text{geom}}$ can be recovered from $wC$. Thus the action of $W_{\text{geom}}$ is strictly transitive on the set of tuples of circles in the packing. 

We will introduce two new spaces related to $C$: the root space and the weight space. The root space contains all Descartes tuples in $\Pa$. The simple roots, which correspond to the dual circles $d_1, \ldots d_n$, also lie in the root space.  The fundamental weights, which correspond to the original circles $c_1, \ldots c_m$, lie in the weight space. Both of these spaces have their own versions of the bilinear form and the Apollonian group.

The root space is $V=(\ker C)^{\perp}$, or, equivalently, the row space of $C$ or column space of $C^T$. Every Descartes tuple can be obtained by multiplying $MC$, where $M$ is the $4 \times 4$ matrix of a M\"obius transformation, and then taking the second row. Thus every Descartes tuple lies in this space. The positive simple roots are $\aalpha_1, \ldots \aalpha_n$ where each $\aalpha_j=-C^T P \mathbf{d}_j$ for a dual circle $\mathbf{d}_j$. Note that the simple roots are not linearly independent (except in the tetrahedral packing case). They are nonzero vectors with all nonnegative entries, so they are positively independent. 

The weight space is $V^* = \R^m/ \ker C$. The fundamental weights are the standard basis vectors $\oomega_1, \ldots \oomega_m$ for $\R^m$. In $V^*$, these are still a spanning set, but they are no longer linearly independent (except in the tetrahedral packing case). The matrix $C$ gives an isomorphism from $V^*$ to $\R^4$, with $\oomega_1, \ldots \oomega_m$ mapping to $\mathbf{c}_1, \ldots \mathbf{c}_m$. The dot product of a Descartes tuple $(b_1, \ldots b_m) \in V$ with a weight $(s_1, \ldots s_m)\in V^*$ is well-defined, and gives a duality between the two spaces.

The bilinear form on the weight space $V^*$ is given by $G=C^T P C$, which is the Gram matrix of the initial configuration of circles. It has entries of $1$ along the diagonal and entries $\leq -1$ off the diagonal. Because $G$ is constructed from $P$, it has three positive, one negative, and $m-4$ zero eigenvalues.

In order to give the bilinear form on the root space, we must choose a right inverse $\tilde{C}$ for $C$. Then the bilinear form on $V$ is given by $\tilde{G}=\tilde{C} P^{-1} \tilde{C}^T$. As a matrix, this depends on the choice of $\tilde{C}$, but it defines a unique bilinear form on $V$. Since $\tilde{C}^T$ is a left inverse for $C^T$, and $V$ is the column space of $C^T$, the vector $\tilde{C}^T \mathbf{b}$ is well-defined for $\mathbf{b} \in V$. Again, $\tilde{G}$ has three positive, one negative, and $m-4$ zero eigenvalues.

The identity $C \tilde{G} C^T = P^{-1}$ implies the following theorem, which is the analogue of the Descartes circle theorem for arbitrary polyhedral packings:
\begin{theorem} \label{descartes}
Every Descartes tuple $\mathbf{b}$ in a polyhedral packing $\Pa$ satisfies 
\begin{equation}\mathbf{b}^T \tilde{G} \mathbf{b} = 0 \end{equation}
\end{theorem}
This is a relation on the tuples of curvatures $b$ that can appear in a polyhedral Descartes configuration. It is possible to give similar relations on the other coordinates $\tilde{b}$, $h_1$, $h_2$ as well. These constitute the analogue of the complex Descartes circle theorem in \cite{MR1903421}. 

The two bilinear forms $G$ and $\tilde{G}$ are duals in the sense that, for $\mathbf{s} \in V^*$ and $\mathbf{b} \in V$, $\mathbf{b} \cdot *$ and $\mathbf{s}^TG$ are equal as linear forms on $V^*$ if and only if $* \cdot \mathbf{s}$ and $\mathbf{b}^T\tilde{G}$ are equal as linear forms on $V$. 

The algebraic Apollonian group acts on the root space; its transpose acts on the weight space. The simple reflections $\sigma_1, \ldots \sigma_n$ on $V$ are reflections along the simple roots. For each $j$, 
\begin{equation} \sigma_j(\bb)=\bb-2(\aalpha_j^T\tilde{G}\bb)\aalpha_j = \bb-2C^TP\mathbf{d}_j\mathbf{d}_j^T\tilde{C}^T\bb \end{equation}
Again, the matrix of $\sigma_j$ depends on the choice of $\tilde{C}$, but $\sigma_j$ is well-defined as a linear map on $V$. The algebraic Apollonian group $W$ is generated by $\sigma_1, \ldots \sigma_n$. We have $\sigma_j = C^T\sigma_{j,\text{geom}}^T \tilde{C}^T$ as mappings on $V$, and so $W=C^T W_{\text{geom}}^T \tilde{C}^T$. This means that $W$ has the same presentation as $W_{\text{geom}}$:
\begin{equation}W = \langle \sigma_1, \ldots \sigma_n | \sigma_j^2=I \rangle \end{equation}
Furthermore, as $W_{\text{geom}}$ preserves the quadratic form $P$, $W$ preserves the quadratic form $\tilde{G}$.

For each $w \in W$, $w^T$ is a mapping on the dual space $V^*$, which preserves the quadratic form $C^T P C$. Each $\sigma_j^T$ is a reflection across the plane perpendicular to the vector $\tilde{C} \mathbf{d}_j \in V^*$. The group $W^T$ is isomorphic to $W_{\text{geom}}$ and the action of $W^T$ on fundamental weights $\oomega_i$ matches the action of $W_{\text{geom}}$ on circles $\mathbf{c}_i$ of the initial tuple. 

We close this section with examples of a basis for $\ker C$, the set of simple roots $\lbrace \aalpha_j \rbrace$ the Gram matrix $G$, and the dual Gram matrix $\tilde{G}$, in the cases of the octahedral and cubic packings. The set of simple reflections $\lbrace \sigma_j \rbrace$ generating the algebraic Apollonian group $W$ can be computed from $\lbrace \aalpha_j \rbrace$ and $\tilde{G}$. Note that these objects do not depend on the choice of initial circles $\mathbf{c}_i$ and dual circles $\mathbf{d}_j$. The matrix $\tilde{G}$ depends on a choice of $\tilde{C}$; we can make a canonical choice using the pseudo-inverse, with the property that $\tilde{C} C$ is the orthogonal projection onto $V$. This removes any dependence on $C$. 

In the case of the octahedral packing, 
\begin{equation*}\ker C = \mathrm{span }( (1, -1, 0, 0, -1, 1)^T, (1, 0, -1, -1, 0, 1)^T),
\end{equation*}
\begin{equation*}\left( \begin{array}{ccc} | & \, & | \\ \aalpha_1 & \cdots & \aalpha_8 \\ | & \, & | \end{array}\right)=2 \sqrt{2}\left(
\begin{array}{cccccccc}
 0 & 0 & 0 & 0 & 1 & 1 & 1 & 1 \\
 0 & 0 & 1 & 1 & 0 & 0 & 1 & 1 \\
 0 & 1 & 0 & 1 & 0 & 1 & 0 & 1 \\
 1 & 0 & 1 & 0 & 1 & 0 & 1 & 0 \\
 1 & 1 & 0 & 0 & 1 & 1 & 0 & 0 \\
 1 & 1 & 1 & 1 & 0 & 0 & 0 & 0 \\
\end{array}
\right) , 
\end{equation*}
\begin{equation*}
G=\left(
\begin{array}{cccccc}
 1 & -1 & -1 & -1 & -1 & -3 \\
 -1 & 1 & -1 & -1 & -3 & -1 \\
 -1 & -1 & 1 & -3 & -1 & -1 \\
 -1 & -1 & -3 & 1 & -1 & -1 \\
 -1 & -3 & -1 & -1 & 1 & -1 \\
 -3 & -1 & -1 & -1 & -1 & 1 \\
\end{array}
\right),
\end{equation*}
\begin{equation*} \tilde{G}=\frac{1}{72}\left(
\begin{array}{cccccc}
 7 & -2 & -2 & -2 & -2 & -11 \\
 -2 & 7 & -2 & -2 & -11 & -2 \\
 -2 & -2 & 7 & -11 & -2 & -2 \\
 -2 & -2 & -11 & 7 & -2 & -2 \\
 -2 & -11 & -2 & -2 & 7 & -2 \\
 -11 & -2 & -2 & -2 & -2 & 7 \\
\end{array}
\right).
\end{equation*}

In the case of the cubic packing, 
\begin{equation*}
\begin{split}\ker C = \mathrm{span }&((1, -1, -1, 1, 0, 0, 0, 0)^T, (1, -1, 0, 0, -1, 1, 0, 0)^T, \\ &(1, 0, -1, 0, -1, 0, 1, 0)^T, (0, 0, 0, 0, 1, -1, -1, 1)^T), \end{split}
\end{equation*}
\begin{equation*}
\left( \begin{array}{ccc} | & \, & | \\ \aalpha_1 & \cdots & \aalpha_6 \\ | & \, & | \end{array}\right)=2 \sqrt{2}\left(
\begin{array}{cccccc}
 0 & 0 & 0 & 1 & 1 & 1 \\
 0 & 0 & 1 & 0 & 1 & 1 \\
 0 & 1 & 0 & 1 & 0 & 1 \\
 0 & 1 & 1 & 0 & 0 & 1 \\
 1 & 0 & 0 & 1 & 1 & 0 \\
 1 & 0 & 1 & 0 & 1 & 0 \\
 1 & 1 & 0 & 1 & 0 & 0 \\
 1 & 1 & 1 & 0 & 0 & 0 \\
\end{array}
\right), \\
\end{equation*}
\begin{equation*}
G=\left(
\begin{array}{cccccccc}
 1 & -1 & -1 & -3 & -1 & -3 & -3 & -5 \\
 -1 & 1 & -3 & -1 & -3 & -1 & -5 & -3 \\
 -1 & -3 & 1 & -1 & -3 & -5 & -1 & -3 \\
 -3 & -1 & -1 & 1 & -5 & -3 & -3 & -1 \\
 -1 & -3 & -3 & -5 & 1 & -1 & -1 & -3 \\
 -3 & -1 & -5 & -3 & -1 & 1 & -3 & -1 \\
 -3 & -5 & -1 & -3 & -1 & -3 & 1 & -1 \\
 -5 & -3 & -3 & -1 & -3 & -1 & -1 & 1 \\
\end{array}
\right),
\end{equation*}
\begin{equation*}
\tilde{G}=\frac{1}{128}\left(
\begin{array}{cccccccc}
 5 & 1 & 1 & -3 & 1 & -3 & -3 & -7 \\
 1 & 5 & -3 & 1 & -3 & 1 & -7 & -3 \\
 1 & -3 & 5 & 1 & -3 & -7 & 1 & -3 \\
 -3 & 1 & 1 & 5 & -7 & -3 & -3 & 1 \\
 1 & -3 & -3 & -7 & 5 & 1 & 1 & -3 \\
 -3 & 1 & -7 & -3 & 1 & 5 & -3 & 1 \\
 -3 & -7 & 1 & -3 & 1 & -3 & 5 & 1 \\
 -7 & -3 & -3 & 1 & -3 & 1 & 1 & 5 \\
\end{array}
\right). \\
\end{equation*}
These computations agree with those of Chait-Roth, Cui, and Stier \cite{MR4017945}, \cite{ChaitRothCuiStier}. The computations for the octahedral packing agree with those of Guettler and Mallows \cite{MR2675919}. Guettler and Mallows give two different quadratic forms satisfied by all sextuples in the octahedral packing; both match our $\tilde{G}$ with different choices of basis for $V$. The computations for the cubic packing agree with those of Stange \cite{MR3814328}, who gives the quadratic form in terms of four of the eight circles of an octuple. 

\section{The Domain of Convergence of $Z(\mathbf{s})$}

Recall that for $\s \in \C^m$, 
\begin{equation} Z(\s) = \sum_{\bb \in \Pa} e^{-\bb \cdot \s} \end{equation}
the sum being over all Descartes $m$-tuples of curvatures $\bb=(b_1, \ldots b_m)$ in a fixed bounded polyhedral packing $\Pa$, with multiplicity. Our goal is to describe the region of absolute convergence of the function $Z(\s)$. This region is a convex tube domain. Since the imaginary part of $\s$ does not affect absolute convergence, we may restrict $\s$ to $\R^m$. Furthermore, using the results of previous section, since all $\bb \in V$, we can view $Z$ as a function on the four-dimensional weight space $V^*$. In the next section, we will visualize the domain of convergence in three-dimensional projective space.

We can rewrite the function $Z(\s)$ using the algebraic Apollonian group. If matrix $C$ represents the initial $m$-tuple of circles in a packing, then $\sigma_j C^T= (\sigma_{j,\text{geom}} C)^T$ for each generator $\sigma_j$ of $W$. Since $W_{\text{geom}}$ acts strictly transitively on the set of $m$-tuples of circles in the packing (see the comment after Lemma \ref{rank4}), $W$ acts strictly transitively on this set as well. Note that the action of $W$ may not be strictly transitive on $m$-tuples of curvatures in the packing; some $m$-tuples of curvatures may appear twice. This will be discussed further in Proposition \ref{multiplicity} below. 

We may write 
\begin{equation} Z(\s)=\sum_{w \in W} e^{-w \bb \cdot \s} =  \sum_{w \in W} e^{- \bb \cdot w^T\s} \end{equation}
where $\bb$ is some fixed $m$-tuple of curvatures in $\Pa$. We see that $Z(\s)$ is symmetric under the action of $W^T$ on $\s$, and the domain of absolute convergence is invariant under this action. 

In order to determine the domain of convergence of $Z(\s)$, we must analyze the action of $W$ on the space of Descartes tuples $V$. Let $\bb \in V$ be an arbitrary Descartes $m$-tuple of curvatures. Define a partial ordering $\leq$ on $V$: $\bb \leq \bb'$ if and only if $\bb'-\bb$ has all nonnegative entries.  For $j=1, \ldots n$, we have $\sigma_j(\bb)=\bb-2(\aalpha_j^T\tilde{G} \bb) \aalpha_j$. Then since $\aalpha_j \geq \mathbf{0}$, we have $\sigma_j(\bb) \geq \bb$ if $\aalpha_j^T\tilde{G} \bb \leq 0$, and $\sigma_j(\bb) \leq \bb$ if $\aalpha_j^T\tilde{G} \bb \geq 0$. The following lemma shows that two distinct generators of $W$ cannot both decrease an $m$-tuple of curvatures in a bounded packing. 

\begin{lemma}\label{decreasing} Suppose that $\aalpha_j^T\tilde{G} \bb \geq 0$ and $\aalpha_k^T\tilde{G} \bb \geq 0$ for distinct positive simple roots $\aalpha_j, \aalpha_k$. Then the sequence 
\begin{equation} \bb, \sigma_j(\bb), \sigma_k\sigma_j(\bb), \sigma_j\sigma_k\sigma_j(\bb), \sigma_k \sigma_j \sigma_k \sigma_j(\bb), \ldots \end{equation}
is monotonically decreasing. In particular, this cannot occur if $\bb$ is an $m$-tuple of curvatures in a bounded polyhedral packing.
\end{lemma}

\begin{proof}
Note that every term in this sequence has the form $\bb+a_j \aalpha_j + a_k \aalpha_k$. We will proceed inductively. Since $\aalpha_j^T\tilde{G} \bb \geq 0$, we have $\sigma_j(\bb) \leq \bb$. Next, suppose that $\bb+a_j \aalpha_j + a_k \aalpha_k$ is an arbitrary term in the sequence, and suppose that the desired inequalities have been verified up to this term. Suppose without loss of generality that $\bb+a_j \aalpha_j + a_k \aalpha_k$ was obtained by applying $\sigma_j$ to the previous term of the sequence. Then $$\sigma_j(\bb+a_j \aalpha_j+a_k \aalpha_k) \geq \bb+a_j \aalpha_j+a_k \aalpha_k$$ so $\aalpha_j^T\tilde{G} \bb + a_j + a_k \aalpha_j^T\tilde{G} \aalpha_k \leq 0$, which means $a_j + a_k \aalpha_j^T\tilde{G} \aalpha_k \leq 0$. Since we have $\aalpha_j^T\tilde{G} \aalpha_k \leq -1$, it follows that $a_j \aalpha_k^T\tilde{G} \aalpha_j + a_k  \geq 0$. This means $\aalpha_k^T\tilde{G} \bb + a_j \aalpha_k^T\tilde{G} \aalpha_j + a_k \geq 0$, and so $$\sigma_k(\bb+a_j \aalpha_j+a_k \aalpha_k) \leq \bb+a_j \aalpha_j+a_k \aalpha_k$$ as desired.  

A monotonically decreasing sequence of $m$-tuples cannot appear in a bounded packing because a bounded packing contains a unique circle of negative curvature, and cannot contain infinitely many different circles with curvature bounded above.  
\end{proof}

Lemma \ref{decreasing} carries useful information about the action of $W$ on $m$-tuples of curvatures in a bounded packing. One example is in the proof of the following proposition. The second part of this proposition is not strictly necessary to the proof of our main theorem, but holds independent interest. It is perhaps surprising that this result holds uniformly for all polyhedra. We give an alternate, more geometric proof of this proposition in Appendix A. 

\begin{prop} \label{multiplicity}
A bounded packing $\Pa$ contains a unique base $m$-tuple of curvatures $\bb$ such that $\aalpha_j^T\tilde{G} \bb \leq 0$ for all $j$. Either every $m$-tuple of curvatures appears with multiplicity 2 in $\Pa$, or every $m$-tuple appears with multiplicity 1. 
\end{prop}
Note that this statement refers to ordered $m$-tuples; unordered tuples can appear with greater multiplicity depending on the symmetry of the polyhedron.
\begin{proof} 
Let $\bb$ be an arbitrary tuple in $\Pa$ and let $\sigma_{j_1}, \sigma_{j_2}, \sigma_{j_3}, \ldots \sigma_{j_k}$ be a sequence of generators for $W$ with consecutive terms distinct. Define $\bb_0=\bb$, $\bb_1=\sigma_{j_1}(\bb)$, $\bb_2=\sigma_{j_2}\sigma_{j_1}(\bb)$, etc. For all $\ell$, we cannot have both $\bb_{\ell} \geq \bb_{\ell-1}$ and $\bb_{\ell} \geq \bb_{\ell+1}$; otherwise $\bb_{\ell}$ would violate Lemma \ref{decreasing}. Therefore we must have $\bb_0>\bb_1>\cdots>\bb_{\ell}$, $\bb_{\ell}\leq\bb_{\ell+1}< \cdots < \bb_k$ for some $\ell$. 

To prove the first statement, choose $\bb \in \Pa$ such that the sum of its entries is minimal. Such a $\bb$ must exist because a bounded packing cannot have an infinite sequence of tuples with the sums of the entries decreasing. Because $\bb \leq \sigma_j(\bb)$, $\aalpha_j^T\tilde{G} \bb \leq 0$ for all $j$. To show uniqueness, suppose that $\bb_0, \bb_1, \bb_2, \ldots \bb_k$ is a sequence as above, and that $\bb_0, \bb_k$ are both base tuples, i.e. $\aalpha_j^T\tilde{G} \bb_0 \leq 0$, $\aalpha_j^T\tilde{G} \bb_k \leq 0$ for all $j$. This is only possible if $k=0$ or $1$, and $\bb_0=\bb_k$.

For the second statement, suppose that a tuple $\bb \in \Pa$ appears with multiplicity greater than 1, and let $\bb_0, \bb_1, \bb_2, \ldots \bb_k$ be a sequence as above with $\bb_0 = \bb_k= \bb$. Then we must have $\sigma_{j_1}=\sigma_{j_k}$, the unique generator which decreases $\bb$. Thus $\bb_1=\bb_{k-1}$, and repeating inductively, $\bb_2=\bb_{k-2}$, etc. We conclude that $k$ must be odd and that the word $\sigma_{j_k} \cdots \sigma_{j_1}$ is a palindrome. Furthermore, this word is uniquely determined, so $\bb$ appears with multiplicity exactly 2 in $\Pa$. Finally, in this case $\bb_{\ell}=\bb_{\ell+1}$; $\sigma_{\ell+1}$ fixes $\bb_{\ell}$ and the other generators increase it, so $\bb_{\ell}$ is the base tuple in $\Pa$ and it appears with multiplicity 2. Since every tuple in $\Pa$ can be obtained from a word in $W$ applied to $\bb_{\ell}$, every tuple appears with multiplicity 2.
\end{proof}

The following proposition gives a lower bound on the growth rate of an $m$-tuple in $\Pa$ as elements of $W$ are applied. 
\begin{prop} \label{growthrate}
Let $\bb \in \Pa$ be the base $m$-tuple of curvatures with $\aalpha_j^T\tilde{G} \bb \leq 0$ for all $j$ and $\aalpha_j^T\tilde{G} \bb = 0$ for at most one $j$. Let $\sigma_{j_1}, \sigma_{j_2}, \sigma_{j_3}, \ldots$ be a sequence of generators for $W$ with consecutive terms distinct.  For $k \in \N$, let 
\begin{equation} \sigma_{j_k}\cdots\sigma_{j_1}(\bb)-\sigma_{j_{k-1}}\cdots\sigma_{j_1}(\bb) = d_k \aalpha_{j_k} \end{equation}
Then $d_k \geq 2 \mu (k-1)$, where $\mu$ is the minimum nonzero value of $|\aalpha_j^T\tilde{G} \bb|$. Thus the total quantity of positive simple roots added to $\bb$ in order to produce $\sigma_{j_k}\cdots\sigma_{j_1}(\bb)$ is at least $\mu k (k-1)$. 
\end{prop}

\begin{proof} 
Let $\sigma_{j_k}\cdots\sigma_{j_1}(\bb) = \bb + a_1 \aalpha_1 + \cdots + a_n \aalpha_n$. Then we have 
\begin{align*}
&d_k = 2\aalpha_{j_k}^T \tilde{G}(\bb + a_1\aalpha_1 + \cdots + a_n \aalpha_n) \\
&d_{k+1} = -2 \aalpha_{j_{k+1}}^T \tilde{G}(\bb + a_1 \aalpha_1 + \cdots + a_n \aalpha_n) \\
&d_{k+1}-d_k = -2(\aalpha_{j_{k+1}}+\aalpha_{j_k})^T \tilde{G}(\bb + a_1 \aalpha_1 + \cdots + a_n \aalpha_n)
\end{align*}
Since $(\aalpha_{j_{k+1}}+\aalpha_{j_k})^T \tilde{G} \aalpha_1, \ldots (\aalpha_{j_{k+1}}+\aalpha_{j_k})^T \tilde{G} \aalpha_n \leq 0$ and $(\aalpha_{j_{k+1}}+\aalpha_{j_k})^T \tilde{G} \bb \leq -\mu$ (and we may assume inductively that $a_1, \ldots a_n \geq 0$), we conclude that $d_{k+1}-d_k \geq 2 \mu$. The proposition follows by induction, since $d_1 \geq 0$. 
\end{proof}

The bound of Proposition \ref{growthrate} is the best possible; indeed, if $\aalpha_1^T \tilde{G} \bb=0$, $\aalpha_2^T \tilde{G} \bb=-\mu$, and $\aalpha_1^T \tilde{G} \aalpha_2=-1$, then for the sequence $\sigma_1, \sigma_2, \sigma_1, \sigma_2, \ldots$, we have $d_k = 2 \mu (k-1)$. This proposition allows us to establish an initial domain of absolute convergence for $Z(\s)$:

\begin{prop}\label{initialdomain}
For a vector $\s$ in the weight space $V^*$, if $\s$ satisfies $\aalpha_j \cdot \s \geq 0$ for all $j$, and $\aalpha_j \cdot \s = 0$ for at most one $j$, then the series $Z$ converges absolutely at $\s$. 
\end{prop}

\begin{proof}
Let $\s \in V^*$ be a weight vector satisfying these hypotheses, and let $\bb \in \Pa$ be the base $m$-tuple, satisfying the hypotheses of Proposition \ref{growthrate}. Then, in the notation of Proposition \ref{growthrate}, 
$$\sigma_{j_k} \cdots \sigma_{j_1}(\bb) \cdot \s = \bb \cdot \s + d_k \aalpha_{j_k} \cdot \s + \cdots + d_1 \aalpha_{j_1} \cdot \s$$
Since at most one $\aalpha_j \cdot \s = 0$, at most every other term in this sum is $0$. Let $\nu$ be the minimum nonzero value of $\aalpha_j \cdot \s$. We conclude by Proposition \ref{growthrate} that 
$$\sigma_{j_k} \cdots \sigma_{j_1}(\bb) \cdot \s \geq \bb \cdot \s + \mu \nu(k-2) +\mu \nu (k-4) + \mu \nu (k-6) + \cdots \geq \bb \cdot \s + \frac{\mu \nu k(k-2)}{4}$$
Now we write $Z (\s)$ as a sum over $w \in W$ and sort by the length $k$ of a reduced word for $w$. There are $s(s-1)^{k-1}$ reduced words of length $k$ for $k \geq 1$. This and the above reasoning give the upper bound 
$$Z(\s) = \sum_{w \in W} e^{-w \bb \cdot \s}  \leq e^{-\bb \cdot \s}\left(1+ \sum_{k=1}^{\infty} s(s-1)^{k-1} e^{-\mu \nu k(k-2)/4} \right)$$
which implies absolute convergence.
\end{proof}

We let $A_0 \subset V^*$ be the domain described in Proposition \ref{initialdomain}, of vectors $\s$ satisfying $\aalpha_j \cdot \s \geq 0$ for all $j$, and $\aalpha_j \cdot \s = 0$ for at most one $j$. We are now able to use the $W^T$-symmetry of $Z(\s)$ to determine its precise domain of absolute convergence:

\begin{theorem}\label{finaldomain}
The domain of absolute convergence of $Z(\s)$ is 
\begin{equation}
A=\bigcup\limits_{w \in W} w^T A_0
\end{equation}
\end{theorem}

\begin{proof}
By Proposition \ref{initialdomain} and $W^T$-invariance, we see that $Z(\s)$ converges absolutely in this region. It suffices to show that $Z(\s)$ diverges elsewhere in $V^*$. 

Suppose that $\s \in V^*$ satisfies $\aalpha_{j_1} \cdot \s, \aalpha_{j_2} \cdot \s \leq 0$, for $j_1 \neq j_2$. Then if $\bb$ is the base tuple of $\Pa$, it follows from Proposition \ref{growthrate} that the sequence $\bb, \sigma_{j_1} \bb, \sigma_{j_2}\sigma_{j_1} \bb, \sigma_{j_1}\sigma_{j_2}\sigma_{j_1}\bb, \sigma_{j_2}\sigma_{j_1}\sigma_{j_2}\sigma_{j_1}\bb$ is monotonically increasing, and the dot products of these vectors with $\s$ are monotonically decreasing. This means that infinitely many terms of the series $Z(\s)$ will be bounded below by $e^{-\bb \cdot \s}$, so this series must diverge.

We have shown that if, for any $w \in W$, $\aalpha_j \cdot w^T \s \geq 0$ for all $j$, and $\aalpha_j \cdot w^T \s = 0$ for at most one $j$, then $Z$ converges at $\s$. Similarly, if for any $w \in W$, $\aalpha_j \cdot w^T \s \leq 0$ for two or more $j$, then $Z$ diverges at $\s$. The only remaining possibility is that for all $w \in W$, $\aalpha_j \cdot w^T \s < 0$ for exactly one $j$ and $\aalpha_j \cdot w^T \s > 0$ for all other $j$. In this case, we can form the unique sequence of generators $\sigma_{j_1}, \sigma_{j_2}, \sigma_{j_3}, \ldots$ such that $\aalpha_{k+1} \cdot \sigma_{j_k}^T \cdots \sigma_{j_1}^T (\s)<0$ for all $k$. We can compute $\aalpha_j \cdot \sigma_j (\s) = -\aalpha_j\cdot \s$, so $\sigma_j$ cannot appear in this sequence twice consecutively. Then for all $k$, we have
\begin{equation*}
\begin{split} \bb \cdot  \sigma_{j_{k+1}}^T \sigma_{j_{k}}^T \cdots \sigma_{j_1}^T (\s)&= \sigma_{j_{k+1}}(\bb) \cdot \sigma_{j_{k}}^T \cdots \sigma_{j_1}^T (\s) \\
&=  (\bb -2 (\aalpha_{k+1}^T \tilde{G} \bb) \aalpha_{k+1}) \cdot \sigma_{j_{k}}^T \cdots \sigma_{j_1}^T (\s) \\
&\leq \bb \cdot  \sigma_{j_{k}}^T \cdots \sigma_{j_1}^T (\s) \end{split}
\end{equation*}
so this sequence of dot products is monotonically decreasing. We again conclude that infinitely many terms of the series $Z(\s)$ will be bounded below by $e^{-\bb \cdot \s}$, so this series must diverge.

\end{proof} 

\section{Geometry of the Domain of Convergence}
In this section, we give a geometric description of the domain of convergence $A$ of Theorem \ref{finaldomain}, generalizing Theorem 5.2 of \cite{Whitehead}. This domain of convergence is the interior of the Tits cone for an infinite root system of complexity beyond affine or hyperbolic types. The upshot of our Theorem \ref{domaingeometry} is that the geometry of the Tits cone reflects the geometry of the underlying packing. Every polyhedral packing type appears in the Tits cone of some root system. 

We argue that the domain of convergence of the generating function is the timelike cone, with infinitely many spacelike protuberances corresponding to the weights in the root system. In an appropriate 3-dimensional projection, this can be visualized as the union of the timelike ball with infinitely many spacelike cones.

We view $A$ as a subset of the four-dimensional weight space $V^*$. This space is equipped with the bilinear form $G$ of signature $(3,1)$, which determines the timelike cone $J=\lbrace \s \in V^* | \s^T G \s < 0 \rbrace$ and its boundary the lightlike cone $N=\lbrace \s \in V^* | \s^T G \s = 0 \rbrace$. Vectors $\s$ satisfying $\s^T G \s >0$ are said to be spacelike. To visualize $A$, we need to further cut the dimension by one. Because $Z(\s)$ is an exponential series, if $\s$ is in the domain of convergence, a multiple of $\s$ by a positive constant will stay in the domain of convergence. Thus, we may visualize $A$ in the projective space $\RP^3$ instead. From the fact that $Z(\s)$ is an exponential series, it also follows that $A$ is convex.

In some basis, the quadratic form defined by $G$ on $V^*$ is equivalent to $x^2+y^2+z^2-t^2$. Then in a projection through an affine plane orthogonal to the $t$ axis, $N$ will appear as a sphere and $J$ as ball. In the rest of this section, we will describe the domains $A_0$ and $A$ as subsets of $\RP^3$. We will implicitly use this basis and projection in our descriptions, e.g. referring to $J$ as the timelike ball, but of course the choice of basis and projection can be changed by applying an automorphism of $\RP^3$.

We will need some geometric properties of the initial domain of convergence $A_0$. The closure $\overline{A_0}$ is a solid convex polyhedron in $\RP^3$, the intersection of the half-spaces $\aalpha_j \cdot \s \geq 0$. The region $A_0$ includes the interior and faces of this polyhedron, but not the edges or the vertices.

\begin{prop} \label{initialgeometry}
The region $\overline{A_0} \subset \RP^3$ is the solid convex polyhedron with vertices at the fundamental weights $\oomega_1, \ldots \oomega_m$. It has the same combinatorial type that defines the packing $\Pa$. All of its edges are tangent to the lightlike sphere $N$. 
\end{prop}

\begin{proof}
From the theory of convex polyhedra, the convex polyhedron with vertices $\oomega_1, \ldots \oomega_m$ is equal to the intersection of the finitely many half-spaces which contain all $\oomega_i$ and have at least three $\oomega_i$ on their boundaries, the faces of the polyhedron. The half-spaces $\aalpha_j \cdot \s \geq 0$ satisfy this condition. Moreover, by equations \eqref{config2}, the vertices $\oomega_i$ and faces $\aalpha_j \cdot \s = 0$ satisfy the correct incidence relations to be a realization of the polyhedron type that defines the packing $\Pa$. Let $\boldsymbol{\beta} \cdot \s \geq 0$ be another half-space which contains all $\oomega_i$, and assume that $\boldsymbol{\beta} \cdot \oomega_i =0$ for some $i$. Because of the incidence relations involving $\oomega_i$ and its neighbors, we have that the half-spaces $\aalpha_j \cdot \s \geq 0$ with $\oomega_i$ on their boundaries form a solid unbounded pyramid with apex $\oomega_i$. The plane  $\boldsymbol\beta \cdot \s = 0$ cannot intersect the interior of this pyramid because it would then separate two vertices. If it intersects the pyramid in just the vertex or in an edge, then it contains only one or two vertices. If it intersects the pyramid in a face $\aalpha_j \cdot \s = 0$, then it coincides with that face. This shows that the inequalities $\aalpha_j \cdot \s \geq 0$ are sufficient to define the polyhedron with vertices $\oomega_1, \ldots \oomega_m$. 

By equations \eqref{config1}, for any two adjacent vertices $\oomega_{i_1}, \oomega_{i_2}$, the quadratic form 
$$(x \oomega_{i_1} + y \oomega_{i_2})^T G (x \oomega_{i_1} + y \oomega_{i_2})=x^2-2xy+y^2$$ is positive semidefinite; this implies that the edge through $\oomega_{i_1}, \oomega_{i_2}$ is tangent to the lightlike sphere. If $\oomega_{i_1}, \oomega_{i_2}$ are not adjacent, then this quadratic form is indefinite and the segment through $\oomega_{i_1}, \oomega_{i_2}$ intersects the timelike ball. 
\end{proof}

We will also need some fundamental properties of the $W^T$ action on $\RP^3$. The generators $\sigma_1^T, \ldots \sigma_n^T$ are projective linear transformations. Because they preserve the bilinear form $G$, they map the lightlike sphere $N$ and timelike ball $J$ to themselves. Finally, each generator $\sigma_j^T$ is a reflection, preserving the plane $\aalpha_j \cdot \s = 0$ and interchanging the two half-spaces cut out by this plane.

The next proposition describes the action of the generators $\sigma_j^T$ restricted to the lightlike sphere $N$, which can be identified with $\hat{\C}$ by stereographic projection. Note that because each circle on $N$ is the intersection of some plane with $N$, and because $\sigma_j^T$ maps planes to planes, it must map circles on $N$ to circles on $N$. 

\begin{prop}\label{reflectiongeometry} Reflections via an algebraic Apollonian group element, $\sigma_j^T$, with corresponding plane of reflection $P: \aalpha_j \cdot \s = 0$, act as circle inversions on $N$, through the circle $c = N \cap P$. \end{prop}

\begin{proof}
It is sufficient to demonstrate two properties: (1) that the two spherical caps on $N$ cut out by $c$ are interchanged, and (2) that all circles on $N$ and perpendicular to $c$ are mapped to themselves.

To see that this is sufficient, consider $T: N \to N$, with these two properties. For any point $\s \in N$, we arbitrarily construct two circles, $d_1$ and $d_2$, perpendicular to $c$ through $\s$. As $T$ preserves circles perpendicular to $c$, we have $T(d_1) = d_1$ and $T(d_2) = d_2$. Note that $d_1$ and $d_2$ will intersect at exactly two points, one of inside $c$ and one outside $c$. As points inside $c$ get mapped to points outside $c$ and $\s$ is one of the intersection points of $d_1$ and $d_2$, $T(\s)$ must be the other. Thus $T(\s)$ is the inversion of $\s$ across $c$.

Property (1) follows from $\sigma_j^T$ being a reflection, fixing $P$ and exchanging the half-spaces cut out by $P$. For property (2), let $d$ be some circle on $N$ perpendicular to $c$, and call either of their points of intersection $\mathbf{t}$. Consider the line $\ell$ tangent to $N$ at $\mathbf{t}$ and also tangent to $c$. As $d$ and $c$ are perpendicular, $\ell$ is also perpendicular to $d$. Now, consider the unique right cone tangent to $N$ with base $d$ and call the apex of this cone $\mathbf{a}$. The segment $\overline{\mathbf{at}}$ is perpendicular to $d$ and by construction, it is also tangent to $N$. Therefore $\overline{\mathbf{at}}$ coincides with $\ell$. Since $\ell$ lies in $P$, $\mathbf{a}$ lies in $P$, and is fixed by $\sigma_j^T$. Finally, as there is a unique cone which is tangent to $N$ with apex $\mathbf{a}$, we have shown that the image of $d$ is $d$ itself, as desired.

\end{proof}

Let $c(\s)$ for spacelike $\s \in \RP^3$ be the circle which is the base of the cone tangent to $N$ with apex $\s$. There is a one-to-one correspondence between spacelike points $\s \in \RP^3$ and circles $c(\s)$ on $N$. Moreover, this correspondence is compatible with the $W^T$-action: since each $w^T \in W^T$ preserves $N$, maps lines to lines, and preserves tangency, it will map the cone with apex $\s$ and base circle $c(\s)$ to the cone with apex $w^T(\s)$ and base circle $w^T(c(\s))$. Then:

\begin{prop} \label{packingonsphere}
The $W^T$-translates of all $c(\oomega_i)$ form a packing $\mathcal{T}$ on the lightlike sphere $N$ of the same polyhedral type as $\mathcal{P}$.
\end{prop}

\begin{proof}
By Proposition \ref{initialgeometry} $\oomega_i$ are the vertices of a polyhedron $\overline{A_0}$ of the type defining $\Pa$ with midsphere $N$. As such, the circles $c(\oomega_i)$ form the initial tuple of a packing of the same type as $\Pa$ on $N$. The Apollonian group generators, $\sigma_i^T$, represent reflections across the faces of the polyhedron $\overline{A_0}$. From Proposition \ref{reflectiongeometry}, we see that those reflections act as circle inversions through the dual configuration to $\lbrace c(\oomega_i) \rbrace$. Having both an initial configuration and inversion in the dual circles, the $W^T$-translates of the circles $c(\oomega_i)$ form a full packing $\mathcal{T}$ of the same type as $\mathcal{P}$ on $N$.
\end{proof}

From this proposition we see that the set of weights $\Omega = \lbrace w^T(\oomega_i) | w \in W, 1 \leq i \leq m \rbrace$, corresponds to the set of circles in $\mathcal{T}$. The weights will be the vertices of the closure of the final domain of convergence $\bar{A}=\overline{\bigcup\limits_{w \in W} w^TA_0}$. We will need the following proposition on the set of weights. The residual set of $\mathcal{T}$ is defined as the sphere $N$ with the disjoint open spherical caps cut out by the circles of $\mathcal{T}$ removed; in particular, it contains all the circles in $\mathcal{T}$. A stronger version of this result is given in \cite[Thm. 3.4]{MR3303040}; we reprove it here in our context.  

\begin{prop}\label{density}
The limit set  of $\Omega$ contains the residual set of $\mathcal{T}$.
\end{prop}

\begin{proof}For any point $\mathbf{t}$ in the residual set, and a given $\epsilon$-neighborhood around it, as the packing is dense on  $N$, there exist infinitely many circles in the packing fully contained in the $\epsilon$-neighborhood. Call one such circle $c$. Assuming that $c$ is small enough that the apex angle of the cone on $c$ is greater than $60^\circ$, then the apex of the cone on $c$ is within distance $\epsilon$ of $\mathbf{t}$. This apex is in $\Omega$, so there exist points of $\Omega$ arbitrarily close to $\mathbf{t}$ 
\end{proof}

\begin{figure}
\includegraphics[width=.72\textwidth]{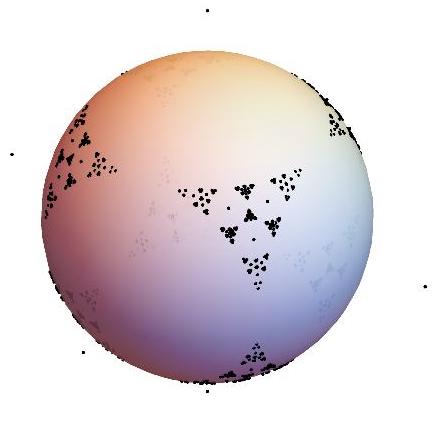}
\caption{The collection of weights for an octahedral packing, shown with the lightlike sphere.}
\label{fig:limitweights}
\end{figure}

The set of weights $\Omega$ for the octahedral packing is illustrated in Figure \ref{fig:limitweights}. We can clearly see that the set of limit weights is the residual set of an octahedral packing on $N$. 

This takes us to the main result of this section:

\begin{theorem} \label{domaingeometry}
For any polyhedral packing, $\mathcal{P}$, the function $Z(\s)$ has a domain of convergence which (under an appropriate choice of basis and projection) can be viewed as the infinite union of the ball $J$, and a set of open, solid cones tangent to $N$ whose bases form a packing of the same type as $\mathcal{P}$ on the surface of $N$.
\end{theorem}

\begin{proof}

For spacelike $\s \in \RP^3$, let $C(\s)$ be the open, solid cone tangent to $N$ with apex $\s$ and base circle $c(\s)$. By Proposition \ref{packingonsphere}, the circles $c(\s)$ for weights $\s \in \Omega$ form a packing of the same type as $\mathcal{P}$ on $N$. Because each $w^T \in W^T$ maps the collection of weights to itself, and correspondingly maps the packing to itself, it must map the collection of cones $C(\s)$ on weights $\s$ to itself. 

We have indicated previously that the weights constitute the zero-skeleton of $\bar{A}$. Because the cones $C(\s)$ are tangent to the sphere, so are the line segments connecting apexes of two cones on adjacent circles. These line segments are on the boundary of the cones, and constitute the one-skeleton of $\bar{A}$. More generally, if $\mathbf{t}_1, \mathbf{t}_2$ are points in two cones $\overline{C(\s_1)}$, $\overline{C(\s_2)}$ with disjoint or tangent base circles, then the line segment through $\mathbf{t}_1, \mathbf{t}_2$ is contained in $\bar{J} \cup \overline{C(\s_1)} \cup \overline{C(\s_2)}$. This shows that the union of $\bar{J}$ with any collection of cones $\overline{C(\s)}$ for $\s \in \Omega$ is convex.

More specifically, consider the initial domain of convergence $A_0$ and the cones $C(\oomega_i)$ for $1 \leq i \leq m$. By Proposition \ref{initialgeometry}, $\overline{A_0}$ is the convex polyhedron with vertices equal to the fundamental weights. Since $\bar{J} \cup \bigcup\limits_{i=1}^{m} \overline{C(\oomega_i)}$ is convex and contains the fundamental weights, we know $\overline{A_0} \subseteq \bar{J} \cup \bigcup\limits_{i=1}^{m} \overline{C(\oomega_i)}$. We have indicated that the initial domain of convergence $A_0$ does not contain the vertices and the edges. A point of $\overline{A_0}$ lies on the boundary of $\bar{J}$ or some $\overline{C(\oomega_i)}$ if and only if it lies on a line segment containing two adjacent $\oomega_i$, i.e. it lies on an edge. Because any point not on the one-skeleton of $\overline{A_0}$ is not on the surface of $\bar{J}$ or any $\overline{C(\oomega_i)}$, we have $A_0\subseteq J\cup \bigcup\limits_{i=1}^{m} C(\oomega_i)$. 

However, the final domain of convergence $A$ is the union of $W^T$ translates of $A_0$, and $J \cup \bigcup\limits_{\s \in \Omega} C(\s)$ is the union of $W^T$ translates of $J \cup \bigcup\limits_{i=1}^{m} C(\oomega_i)$. Therefore, we have $A \subseteq J\cup \bigcup\limits_{\s \in \Omega} C_s$, as desired. 

For the opposite inclusion, we must show that $C(\s) \subseteq A$ for all $\s \in \Omega$ and that $J \subseteq A$. To show that $C(\s) \subseteq A$, note that the apex $\s$ is a weight, and is on the boundary of $A$. The circle $c(\s)$ is part of the residual set of the packing $\mathcal{T}$, and therefore by Proposition \ref{density}, it also lies in the closure of $\Omega$, and hence of $A$. The convex hull of $\s$ and $c(\s)$ is the cone $\overline{C(\s)}$. Then because the closure of $A$ is also convex, we have $\overline{C(\s)} \subseteq \bar{A}$. Finally, since $C(\s)$ is open and convex and $A$ is convex, $C_s \subseteq A$.

To show that $J \subseteq A$, notice that the open spherical cap cut out of $N$ by each circle in the packing $\mathcal{T}$ lies in $A$ because of the previous paragraph. Because the packing is dense, these spherical caps are dense on $N$. Thus $N$ is contained in $\bar{A}$, and by convexity, $J \subseteq \bar{A}$. Again, since $J$ is open and convex and $A$ is convex, $J \subseteq A$.
\end{proof}

\begin{figure}
\includegraphics[width=.36\textwidth]{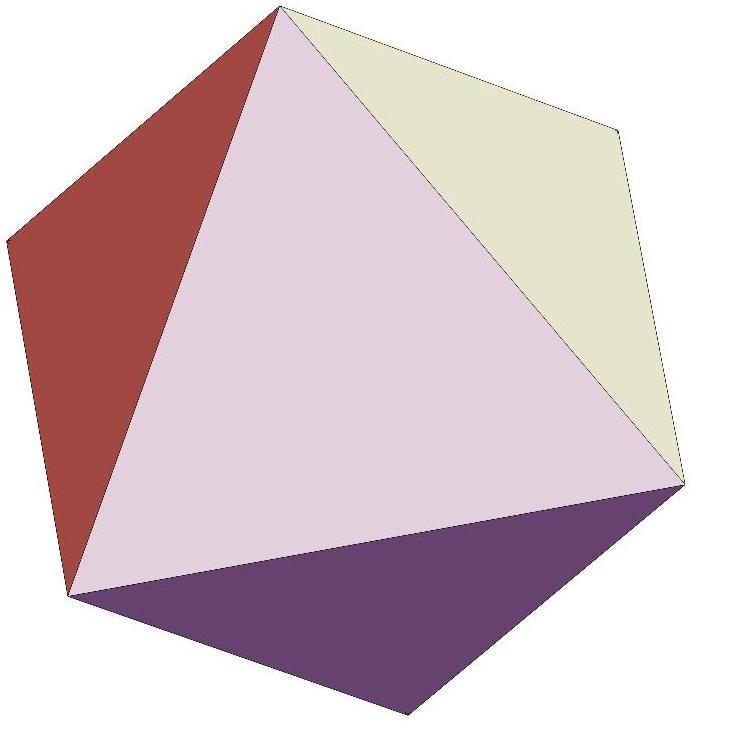} \hspace{.1\textwidth} \includegraphics[width=.36\textwidth]{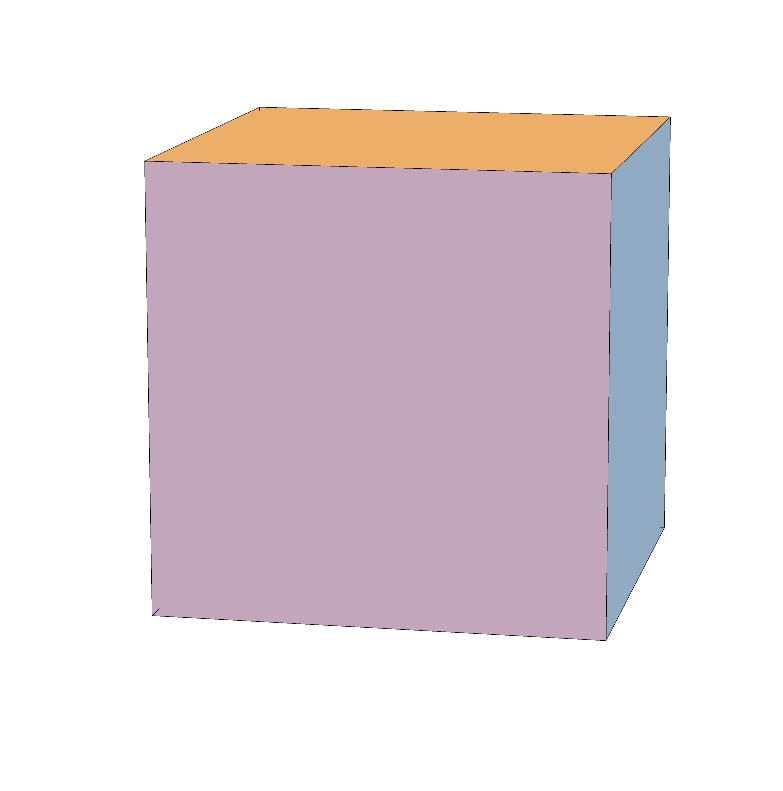}

\includegraphics[width=.36\textwidth]{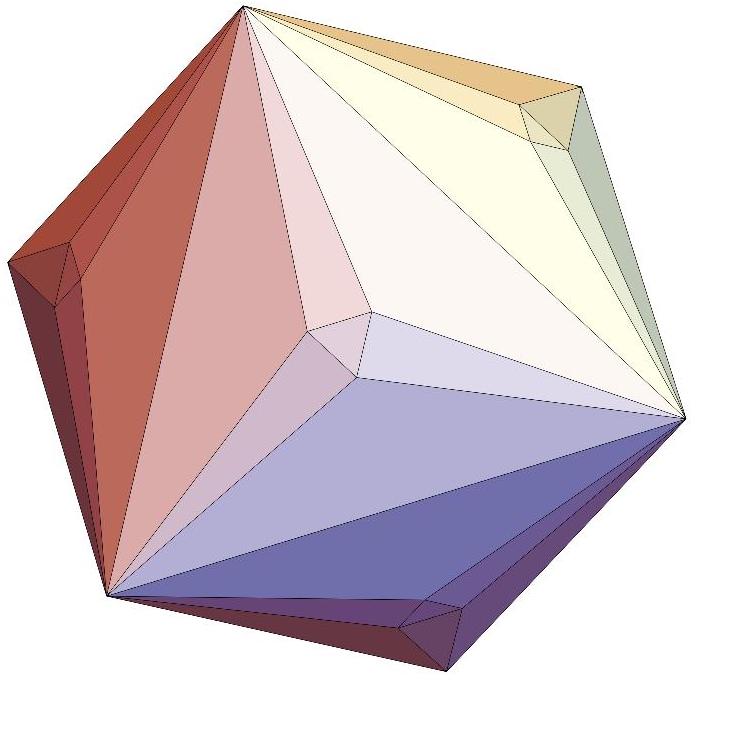} \hspace{.1\textwidth} \includegraphics[width=.36\textwidth]{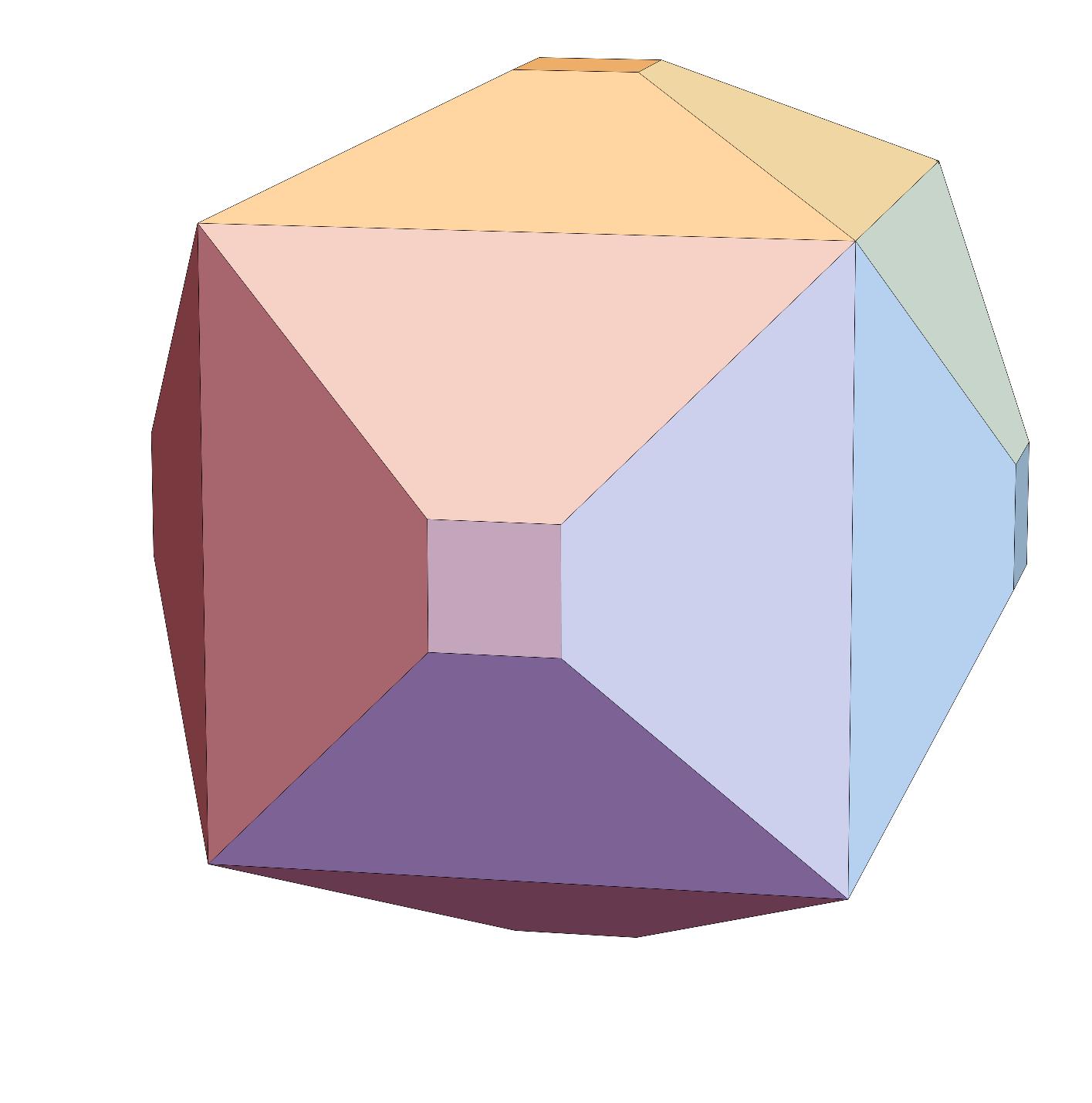}

\includegraphics[width=.36\textwidth]{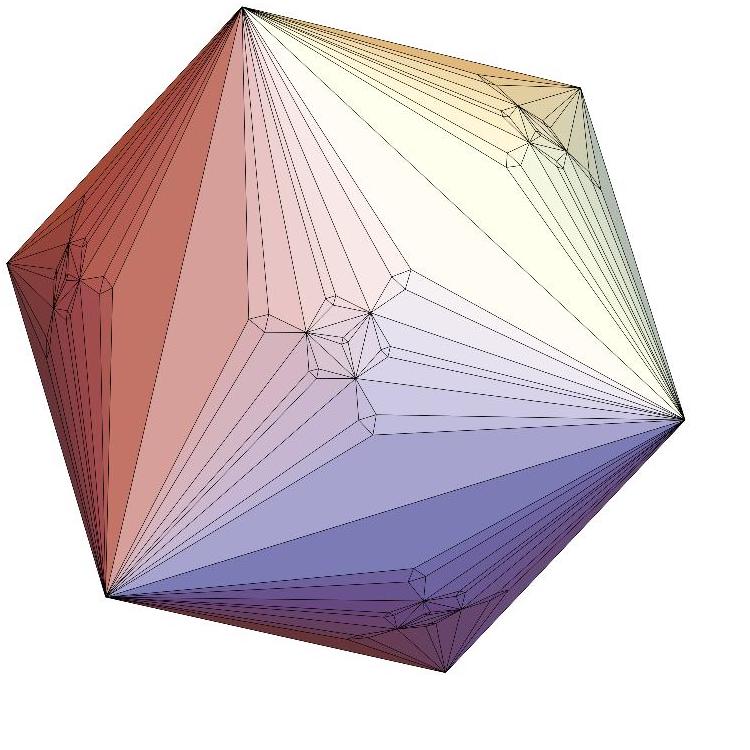} \hspace{.1\textwidth} \includegraphics[width=.36\textwidth]{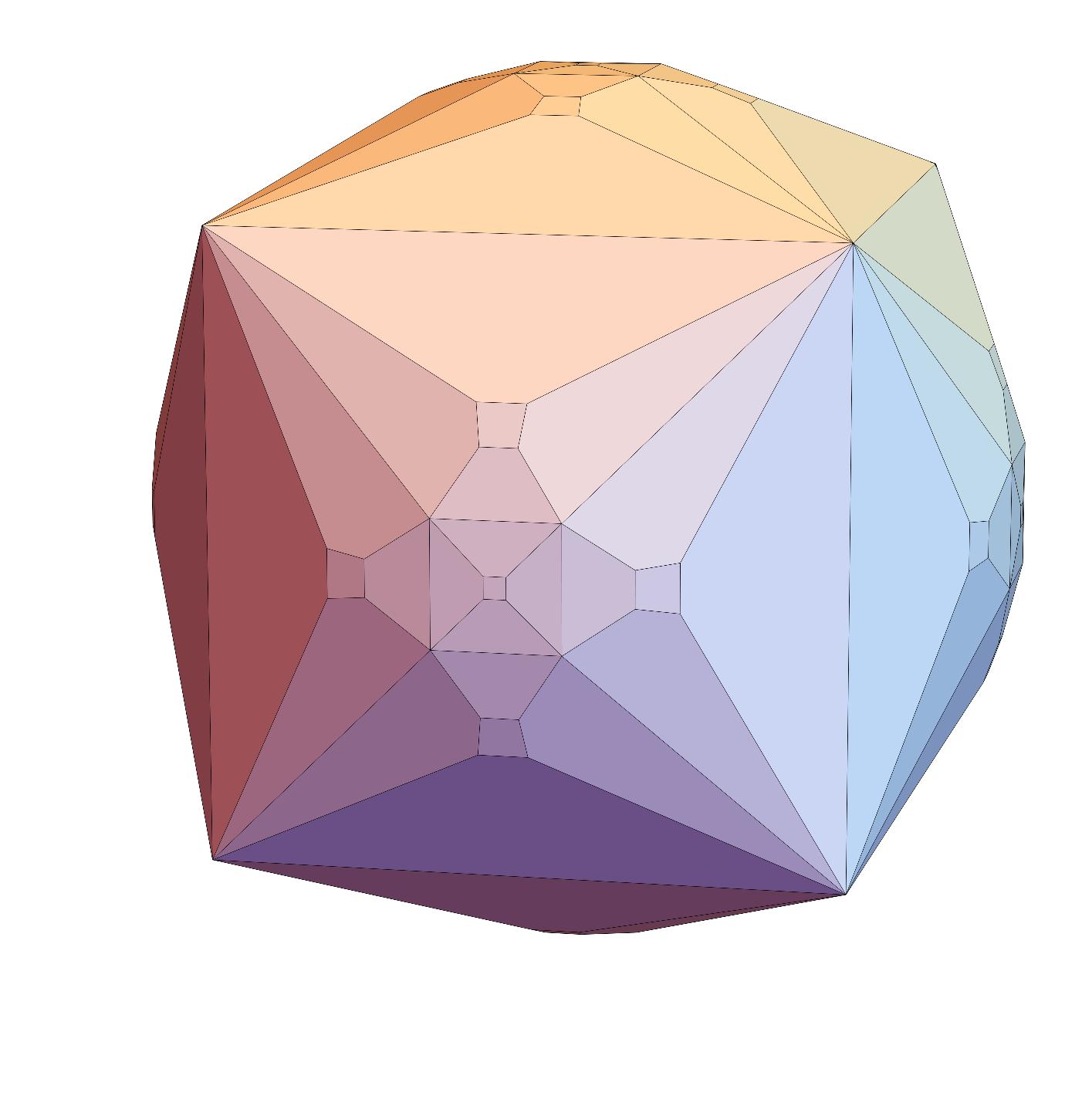}

\includegraphics[width=.36\textwidth]{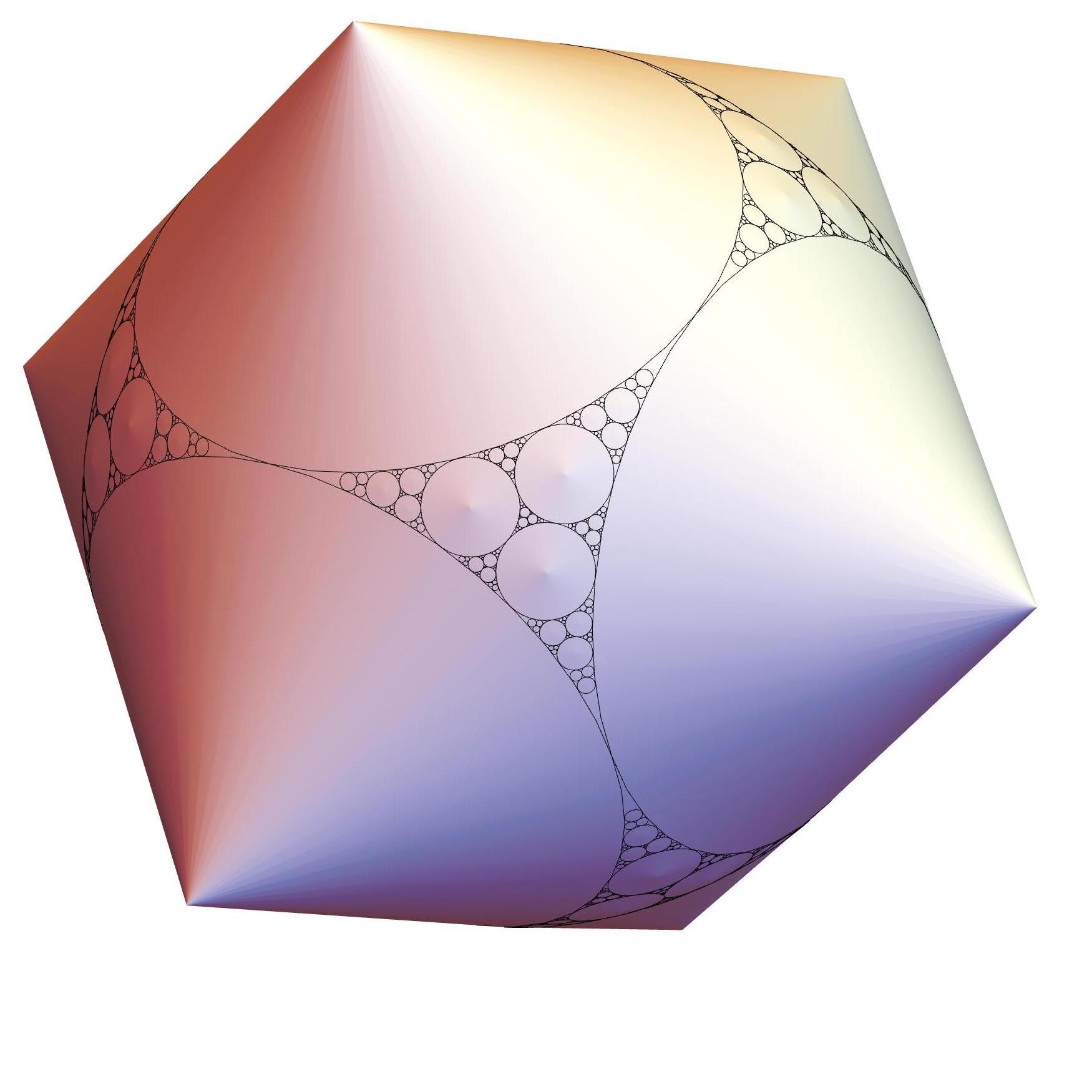} \hspace{.1\textwidth} \includegraphics[width=.36\textwidth]{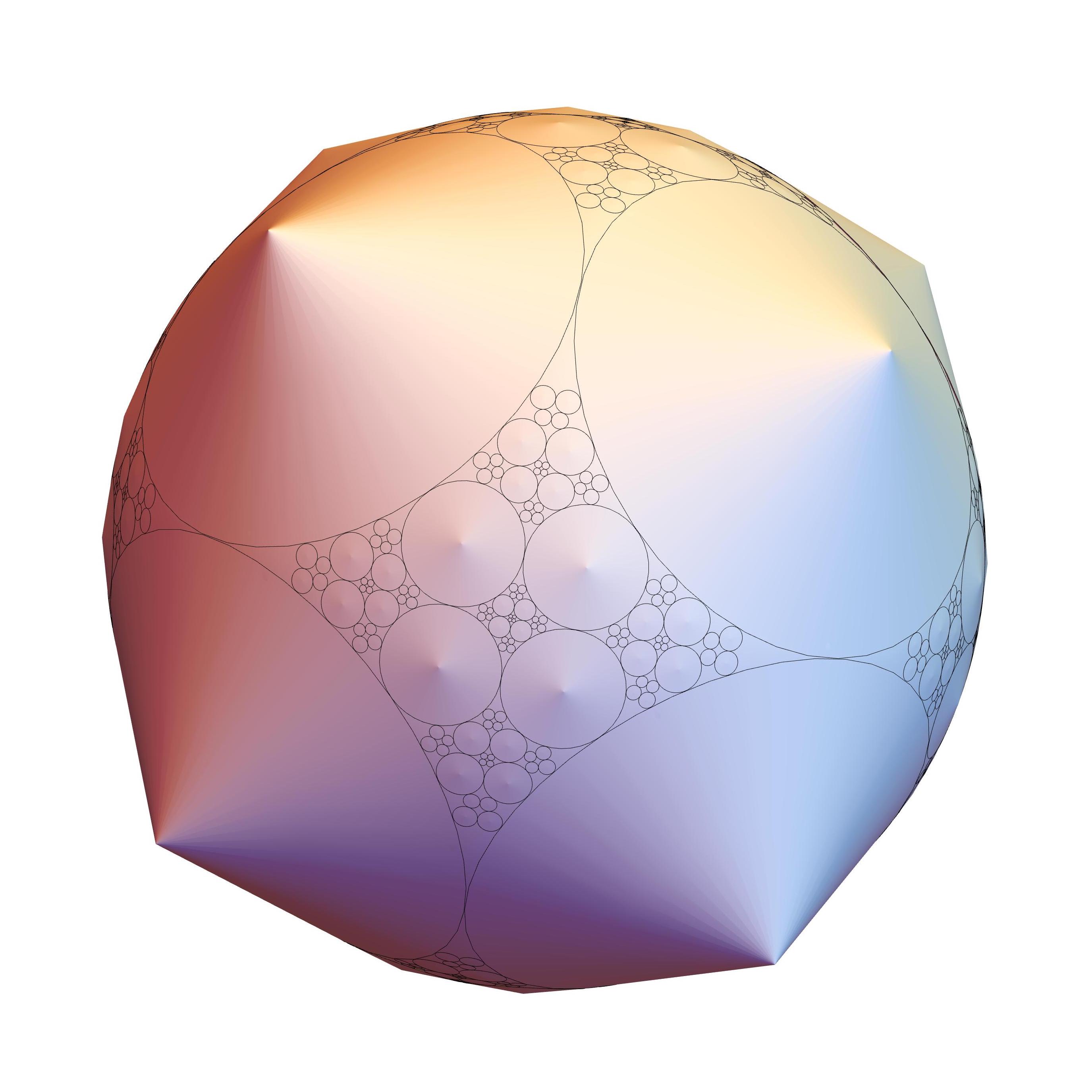}
\caption{The initial domain of convergence, its orbit under words of length 1 and 2 in the Apollonian group, and the final domain of convergence for the octahedral and cubic packings.}
\label{fig:domains}
\end{figure}

Figure \ref{fig:domains} shows the initial polyhedron of convergence $A_0$, its orbit under two generations of generators $\sigma_j^T$ for $W^T$, and the final domain of convergence $A$ for the octahedral and cubic packings. These are two-dimensional pictures of a three-dimensional projection of a four-dimensional quotient of an $m$-dimensional real part of a $2m$-dimensional complex domain. 

\appendix 

\section{A Geometric Proof of Proposition \ref{multiplicity}}

In this appendix, we sketch the proof of a slight modification of Proposition \ref{multiplicity}, using the geometry of polyhedral tuples and dual tuples instead of the algebra of the root space. The proof relies on the following important proposition:

\begin{prop}\label{rigidity}
Let $\bb=(b_1,\ldots b_m)$ be a Descartes $m$-tuple of curvatures which can be realized as a configuration of circles in the plane with the correct tangency relationships and dual circles determined by a polyhedron $\Pi$. Then this realization is unique up to rigid motions. 
\end{prop}

The idea of the proof is to work with one face of the polyhedron at a time. If the curvatures of the circles around a face are all known, then the curvature of the dual circle corresponding to that face can be computed, explicitly and uniquely. Once the size of the dual circle is known, the placement of the ring of circles orthogonal to it is uniquely determined up to rigid motions. Two such rings of circles corresponding to adjacent faces share two tangent circles in common, so they can be glued together in a unique way up to rigid motions and a choice of orientation. Further, only one choice of orientation allows the dual circles to satisfy the correct tangencies. Iteratively gluing together the rings of circles will give the full configuration up to rigid motions, as every circle is orthogonal to some dual circle. 

Our method of proof does not establish that every $\bb \in V$ satisfying $\bb^T \tilde{G} \bb = 0$ can be realized as a polyhedral circle configuration, though we expect this to be true. 

Because the initial configuration and dual configuration are unique up to rigid motions, the entire packing $\Pa$ containing $\bb$ is uniquely determined up to rigid motions. In particular, if $\bb$ is part of a bounded packing, then the size of the external circle for this packing is determined by $\bb$.

\begin{figure}
\includegraphics[width=.36\textwidth]{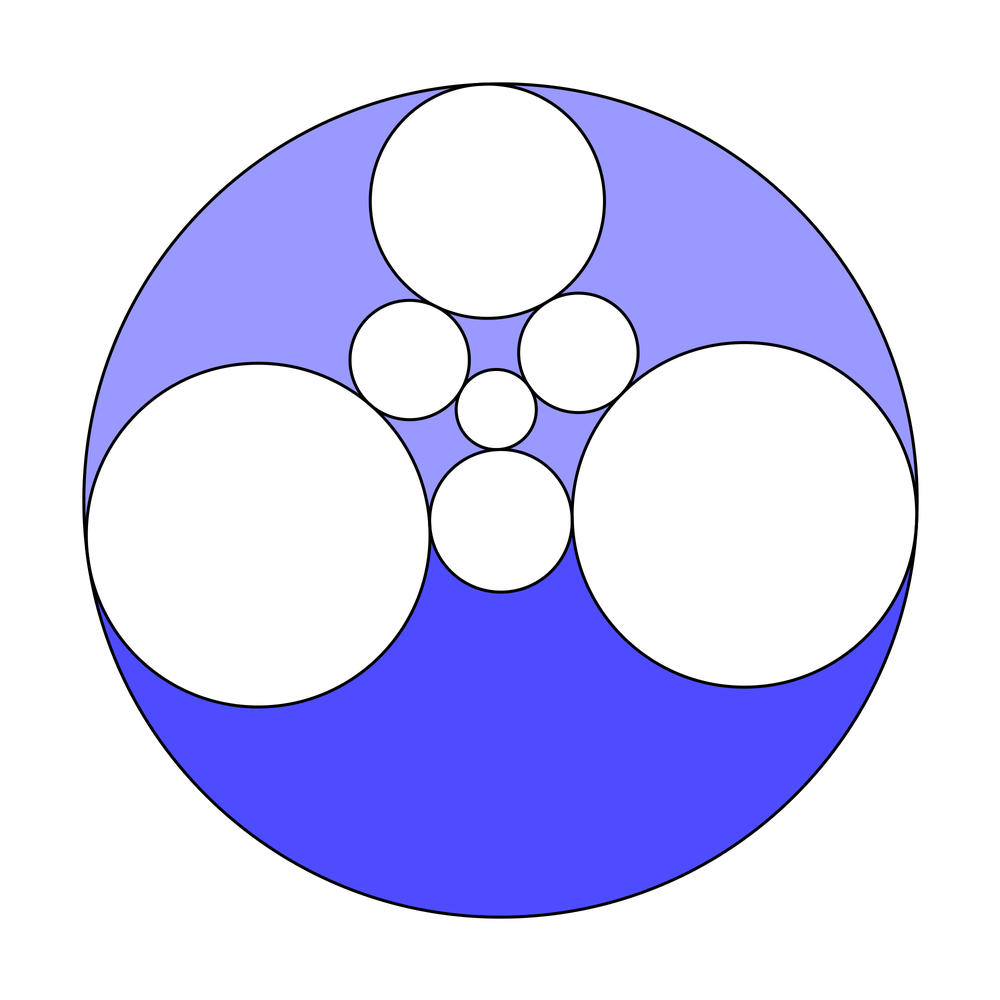} \hspace{.1\textwidth} \includegraphics[width=.36\textwidth]{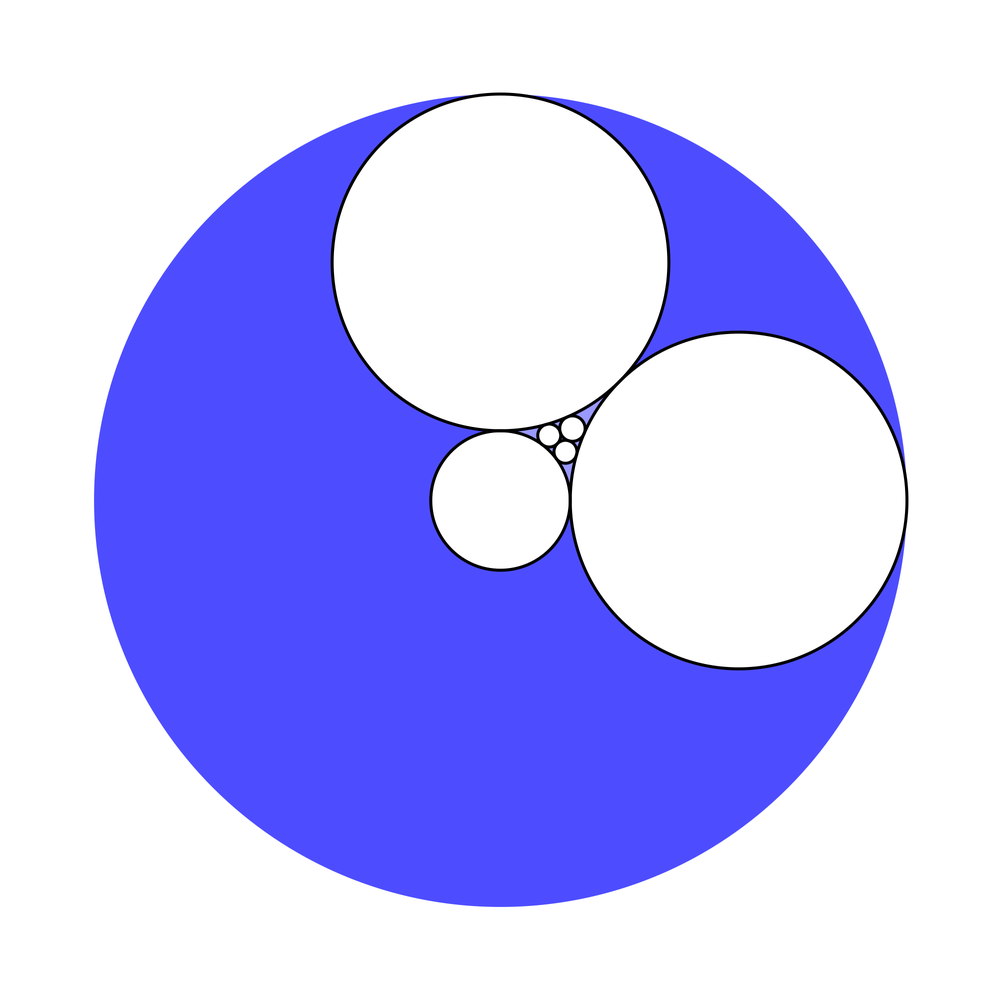}
\caption{The area invariant for a tuple in a cubic packing and a tuple in an octahedral packing. Each region $F_j$ is shaded blue and the largest region is shaded darker blue. In the second image, the largest region is bounded by three circles from the tuple and the exterior circle, and it has two connected components.} 
\label{fig:area}
\end{figure}

We now define an area invariant associated to any tuple $\bb$ in a bounded packing. By Proposition \ref{rigidity}, $\bb$ can be realized with a tuple of circles uniquely up to rigid motion. The circles can be oriented so that their interiors are disjoint. Consider the complement of these circles and their interiors in the plane $\C$. This open region has connected components $F_1, \ldots F_n$, one for each face of the polyhedron $\Pi$. Each connected component $F_j$ is bounded by arcs along the ring of circles for face $j$. The intersection points of these arcs lie on the dual circle $d_j$ corresponding to face $j$, which encloses $F_j$. If $\bb$ contains an external circle, then each region $F_j$ is bounded, but if $\bb$ does not contain an external circle, then exactly one region $F_j$ is unbounded, and contains the point at infinity in $\hat{\C}$. In the latter case, let $E$ be the open disc whose bounding circle is the external circle for the packing that contains $\bb$. According to Proposition \ref{rigidity}, the size of $E$ and its placement relative to the tuple $\bb$ is uniquely determined. We can then replace each $F_j$ by its intersection with $E$, even though this may disconnect one region $F_j$ The use of the external circle is necessary to give each region $F_j$ a finite area. We define
\begin{equation} \mathrm{Area}(\bb) = \max_{1 \leq j \leq n} \mathrm{Area}(F_j)\end{equation}
Proposition \ref{rigidity} ensures that this value is well-defined. Figure \ref{fig:area} illustrates the area invariant for two examples of Descartes tuples, one from a cubic and one from an octahedral packing. 

The following proposition relates the area invariant to the action of $W$: 

\begin{prop}\label{onegenerator}
There is at most one generator $\sigma_j$ of $W$ such that $\mathrm{Area}(\sigma_j (\bb)) \leq \mathrm{Area}(\bb)$.
\end{prop}

The generator $\sigma_j$ can be realized as an inversion through the dual circle $d_j$ corresponding to face $j$. This inversion maps all connected components $F_k$ for $k \neq j$ to subsets of $F_j$; it maps $F_j$ to a set containing all the other components $F_k$. This means that in order for $\sigma_j$ to decrease $\mathrm{Area}(\bb)$, the area of $F_j$ must be larger than the sum of the areas of all the other connected components. This can hold for at most one $j$.

With Proposition \ref{onegenerator} established, we can order tuples $\bb \in \Pa$ by area rather than by the partial ordering of Section 3. If $\bb_0=\bb$, $\bb_1=\sigma_{j_1}(\bb)$, $\bb_2=\sigma_{j_2}\sigma_{j_1}(\bb)$, etc., then we must have $\mathrm{Area}(\bb_0)>\mathrm{Area}(\bb_1)>\cdots>\mathrm{Area}(\bb_{\ell})$, $\mathrm{Area}(\bb_{\ell})\leq \mathrm{Area}(\bb_{\ell+1})< \cdots < \mathrm{Area}(\bb_k)$ for some $\ell$; otherwise, some tuple has its area reduced by two distinct generators. The proof of the second statement of Proposition \ref{multiplicity} carries through exactly as before. 

The first statement requires some modification. Let us define a geometric base tuple in $\Pa$ as a tuple with minimal area. This does not necessarily coincide with the base tuple defined in Proposition \ref{multiplicity}. We will show that a bounded polyhedral packing $\Pa$ contains either one or two geometric base tuples. 

If $\Pa$ contains no geometric base tuple, then there is an infinite sequence of tuples $\bb_0=\bb, \, \bb_1=\sigma_{j_1}(\bb), \, \bb_2=\sigma_{j_2}\sigma_{j_1}(\bb), \ldots$ with decreasing area. As in the proof of Proposition \ref{onegenerator}, each tuple $\bb_k$ has a region $F_{j_{k+1}}(\bb_k)$ whose area is greater than the areas of all the other regions combined. The inversion $\sigma_{j_{k+1}}$ maps the other regions to subsets of $F_{j_{k+1}}(\bb_k)$ and maps $F_{j_{k+1}}(\bb_k)$ to a superset of the other regions. Since $\sigma_{j_{k+2}} \neq \sigma_{j_{k+1}}$, the maximal region at the next step, $F_{j_{k+2}}(\bb_{k+1})$, must be a subset of $F_{j_{k+1}}(\bb_k)$. Thus the maximal regions $F_{j_{k+1}}(\bb_k)$ are nested, with area bounded below by a positive constant. They are contained in nested dual circles $d_{j_{k+1}}(\bb_k)$, with radii bounded below by a positive constant. Choose any circle in $\Pa$ outside the outermost dual circle $d_{j_1}(\bb_0)$. Its reflections across the nested dual circles $d_{j_{k+1}}(\bb_k)$ are an infinite collection of circles in $\Pa$ with radius bounded below, a contradiction since $\Pa$ is bounded.

Having shown that a geometric base tuple exists, the proof that there are at most two geometric base tuples carries through as in the proof of Proposition \ref{multiplicity}. It can happen that two distinct tuples, related by an Apollonian group generator $\sigma_j$, have equal, minimal areas.

\bibliography{ApollonianBib}
\bibliographystyle{amsplain}

\end{document}